\documentclass[11pt]{article}
\usepackage{pifont}
\usepackage{syntonly}
\usepackage{mathrsfs}

\usepackage[centertags]{amsmath}
\usepackage{amsfonts}
\usepackage{amssymb}
\usepackage{amsthm}
\usepackage{newlfont}
\usepackage{CJK}
\usepackage{indentfirst}
    \usepackage[colorlinks,linkcolor=blue, pdfstartview=FitH]{hyperref}
\usepackage{url}
\usepackage{color,titlesec,titletoc}
\newlength{\defbaselineskip}
\newcommand{\setlinespacing}[1]%
           {\setlength{\baselineskip}{#1 \defbaselineskip}}

 \theoremstyle{plain}
\newtheorem{thm}{Theorem}[section]
\newtheorem{theorem}[thm]{Theorem}
\newtheorem{lemma}[thm]{Lemma}
\newtheorem{corollary}[thm]{Corollary}
\newtheorem{proposition}[thm]{Proposition}

\theoremstyle{definition}

\newtheorem{remark}[thm]{Remark}

\newtheorem{definition}[thm]{Definition}

\newtheorem{example}[thm]{Example}

\newtheorem{defn-thm}[thm]{Definition-Theorem}

\numberwithin{equation}{section}
\catcode`\@=11
\def\opn#1#2{\def#1{\mathop{\kern0pt\fam0#2}\nolimits}}
\def\underrightarrow{\mathpalette\underrightarrow@}
\def\underrightarrow@#1#2{\vtop{\ialign{$##$\cr
 \hfil#1#2\hfil\cr\noalign{\nointerlineskip}%
 #1{-}\mkern-6mu\cleaders\hbox{$#1\mkern-2mu{-}\mkern-2mu$}\hfill
 \mkern-6mu{\to}\cr}}}

\def\underleftarrow{\mathpalette\underleftarrow@}
\def\underleftarrow@#1#2{\vtop{\ialign{$##$\cr
 \hfil#1#2\hfil\cr\noalign{\nointerlineskip}#1{\leftarrow}\mkern-6mu
 \cleaders\hbox{$#1\mkern-2mu{-}\mkern-2mu$}\hfill
 \mkern-6mu{-}\cr}}}
\let\amp@rs@nd@\relax
\newdimen\ex@
\newdimen\bigaw@
\newdimen\minaw@
\newdimen\minCDaw@
\newif\ifCD@
\def\minCDarrowwidth#1{\minCDaw@#1}

\def\@CD{\def\A##1A##2A{\llap{$\vcenter{\hbox
 {$\scriptstyle##1$}}$}\Big\uparrow\rlap{$\vcenter{\hbox{%
$\scriptstyle##2$}}$}&&}%
\def\V##1V##2V{\llap{$\vcenter{\hbox
 {$\scriptstyle##1$}}$}\Big\downarrow\rlap{$\vcenter{\hbox{%
$\scriptstyle##2$}}$}&&}%
\def\={&\hskip.5em\mathrel
 {\vbox{\hrule width\minCDaw@\vskip3\ex@\hrule width
 \minCDaw@}}\hskip.5em&}%
\def\verteq{\Big\Vert&&}%
\def\noarr{&&}%
\def\vspace##1{\noalign{\vskip##1\relax}}\relax\let\amp@rs@nd@&\iffalse}\fi
 \CD@true\vcenter\bgroup\relax\let\\=\cr\iffalse}\fi\tabskip\z@skip\baselineskip20\ex@
 &\hfill$\m@th##$\hfill\cr}
\def\@endCD{\cr\egroup\egroup}
\def\>#1>#2>{\amp@rs@nd@\setbox\z@\hbox{$\scriptstyle
 \;{#1}\;\;$}\setbox\@ne\hbox{$\scriptstyle\;{#2}\;\;$}\setbox\tw@
 \hbox{$#2$}\ifCD@
 \global\bigaw@\minCDaw@\else\global\bigaw@\minaw@\fi
 \ifdim\wd\z@>\bigaw@\global\bigaw@\wd\z@\fi
 \ifdim\wd\@ne>\bigaw@\global\bigaw@\wd\@ne\fi
 \ifCD@\hskip.5em\fi
 \ifdim\wd\tw@>\z@
 \mathrel{\mathop{\hbox to\bigaw@{\rightarrowfill}}\limits^{#1}_{#2}}\else
 \mathrel{\mathop{\hbox to\bigaw@{\rightarrowfill}}\limits^{#1}}\fi
 \ifCD@\hskip.5em\fi\amp@rs@nd@}
\def\<#1<#2<{\amp@rs@nd@\setbox\z@\hbox{$\scriptstyle
 \;\;{#1}\;$}\setbox\@ne\hbox{$\scriptstyle\;\;{#2}\;$}\setbox\tw@
 \hbox{$#2$}\ifCD@
 \global\bigaw@\minCDaw@\else\global\bigaw@\minaw@\fi
 \ifdim\wd\z@>\bigaw@\global\bigaw@\wd\z@\fi
 \ifdim\wd\@ne>\bigaw@\global\bigaw@\wd\@ne\fi
 \ifCD@\hskip.5em\fi
 \ifdim\wd\tw@>\z@
 \mathrel{\mathop{\hbox to\bigaw@{\leftarrowfill}}\limits^{#1}_{#2}}\else
 \mathrel{\mathop{\hbox to\bigaw@{\leftarrowfill}}\limits^{#1}}\fi
 \ifCD@\hskip.5em\fi\amp@rs@nd@}

\def\@CDS{\def\A##1A##2A{\llap{$\vcenter{\hbox
 {$\scriptstyle##1$}}$}\Big\uparrow\rlap{$\vcenter{\hbox{%
$\scriptstyle##2$}}$}&}%
\def\V##1V##2V{\llap{$\vcenter{\hbox
 {$\scriptstyle##1$}}$}\Big\downarrow\rlap{$\vcenter{\hbox{%
$\scriptstyle##2$}}$}&}%
\def\={&\hskip.5em\mathrel
 {\vbox{\hrule width\minCDaw@\vskip3\ex@\hrule width
 \minCDaw@}}\hskip.5em&}
\def\verteq{\Big\Vert&}
\def\novarr{&}
\def\noharr{&&}
\def\SE##1E##2E{\slantedarrow(0,18)(4,-3){##1}{##2}&}
\def\SW##1W##2W{\slantedarrow(24,18)(-4,-3){##1}{##2}&}
\def\NE##1E##2E{\slantedarrow(0,0)(4,3){##1}{##2}&}
\def\NW##1W##2W{\slantedarrow(24,0)(-4,3){##1}{##2}&}
\def\slantedarrow(##1)(##2)##3##4{%
\thinlines\unitlength1pt\lower 6.5pt\hbox{\begin{picture}(24,18)%
\put(##1){\vector(##2){24}}%
\put(0,8){$\scriptstyle##3$}%
\put(20,8){$\scriptstyle##4$}%
\end{picture}}}
\def\vspace##1{\noalign{\vskip##1\relax}}\relax\let\amp@rs@nd@&\iffalse}\fi
 \CD@true\vcenter\bgroup\relax\let\\=\cr\iffalse}\fi\tabskip\z@skip\baselineskip20\ex@
 &\hfill$\m@th##$\hfill\cr}
\def\@endCDS{\cr\egroup\egroup}
\newdimen\TriCDarrw@
\newif\ifTriV@

\def\@TriCDV{\TriV@true\def\TriCDpos@{6}\@TriCD}
\def\@TriCDA{\TriV@false\def\TriCDpos@{10}\@TriCD}
\def\@TriCD#1#2#3#4#5#6{%
\setbox0\hbox{$\ifTriV@#6\else#1\fi$} \TriCDarrw@=\wd0
\advance\TriCDarrw@ 24pt \advance\TriCDarrw@ -1em
\def\SE##1E##2E{\slantedarrow(0,18)(2,-3){##1}{##2}&}
\def\SW##1W##2W{\slantedarrow(12,18)(-2,-3){##1}{##2}&}
\def\NE##1E##2E{\slantedarrow(0,0)(2,3){##1}{##2}&}
\def\NW##1W##2W{\slantedarrow(12,0)(-2,3){##1}{##2}&}
\def\slantedarrow(##1)(##2)##3##4{\thinlines\unitlength1pt
\lower 6.5pt\hbox{\begin{picture}(12,18)%
\put(##1){\vector(##2){12}}%
\put(-4,\TriCDpos@){$\scriptstyle##3$}%
\put(12,\TriCDpos@){$\scriptstyle##4$}%
\end{picture}}}
\def\={\mathrel {\vbox{\hrule
   width\TriCDarrw@\vskip3\ex@\hrule width
   \TriCDarrw@}}}
\def\>##1>>{\setbox\z@\hbox{$\scriptstyle
 \;{##1}\;\;$}\global\bigaw@\TriCDarrw@
 \ifdim\wd\z@>\bigaw@\global\bigaw@\wd\z@\fi
 \hskip.5em
 \mathrel{\mathop{\hbox to \TriCDarrw@
{\rightarrowfill}}\limits^{##1}}
 \hskip.5em}
\def\<##1<<{\setbox\z@\hbox{$\scriptstyle
 \;{##1}\;\;$}\global\bigaw@\TriCDarrw@
 \ifdim\wd\z@>\bigaw@\global\bigaw@\wd\z@\fi
 \mathrel{\mathop{\hbox to\bigaw@{\leftarrowfill}}\limits^{##1}}
 }
 \CD@true\vcenter\bgroup\relax\let\\=\cr\iffalse}\fi
 \tabskip\z@skip\baselineskip20\ex@
 \ifTriV@
 \halign\bgroup
 &\hfill$\m@th##$\hfill\cr
#1&\multispan3\hfill$#2$\hfill&#3\\
&#4&#5\\
&&#6\cr\egroup%
\else
 \halign\bgroup
 &\hfill$\m@th##$\hfill\cr
&&#1\\%
&#2&#3\\
#4&\multispan3\hfill$#5$\hfill&#6\cr\egroup \fi}
\def\@endTriCD{\egroup}

\newcounter{Myenumi}
{\begin{list}{}{\usecounter{Myenumi}%
\settowidth{\leftmargin}{2.n}\settowidth{\labelwidth}{2.n}%
\setlength{\labelsep}{0pt}}}{\end{list}}
\newcounter{Myenumii}
{\begin{list}{}{\usecounter{Myenumii}%
\settowidth{\leftmargin}{a)n}\settowidth{\labelwidth}{a)n}%
\setlength{\labelsep}{0pt}}}{\end{list}}
\newcounter{Myenumiii}
{\begin{list}{}{\usecounter{Myenumiii}%
\settowidth{\leftmargin}{iv.n}\settowidth{\labelwidth}{iv.n}%
\setlength{\labelsep}{0pt}}}{\end{list}}

{\end{list}}

{\begin{list}{}{%
\settowidth{\leftmargin}{2.n}\settowidth{\labelwidth}{2.n}%
\setlength{\labelsep}{0pt}}}{\end{list}}

\newsymbol\onto 1310

\begin{document}

\title{Geometry of  Hermitian manifolds}
\author{Kefeng Liu,  Xiaokui Yang
}

\date{}

\maketitle
\begin{abstract} On Hermitian manifolds, the second Ricci curvature tensors
of various metric connections are closely related to the geometry of
 Hermitian manifolds.  By refining the Bochner formulas for any
Hermitian complex vector bundle (Riemannain real vector bundle) with
an arbitrary metric connection over a compact Hermitian manifold, we
can derive various vanishing theorems for Hermitian manifolds and
complex vector bundles by the second Ricci curvature tensors. We
will also introduce a natural geometric flow on Hermitian manifolds
by using the second Ricci curvature tensor.
\end{abstract}


\newcommand{\sA}{{\mathcal A}}
\newcommand{\sB}{{\mathcal B}}
\newcommand{\sC}{{\mathcal C}}
\newcommand{\sD}{{\mathcal D}}
\newcommand{\sE}{{\mathcal E}}
\newcommand{\sF}{{\mathcal F}}
\newcommand{\sG}{{\mathcal G}}
\newcommand{\sH}{{\mathcal H}}
\newcommand{\sI}{{\mathcal I}}
\newcommand{\sJ}{{\mathcal J}}
\newcommand{\sK}{{\mathcal K}}
\newcommand{\sL}{{\mathcal L}}
\newcommand{\sM}{{\mathcal M}}
\newcommand{\sN}{{\mathcal N}}
\newcommand{\sO}{{\mathcal O}}
\newcommand{\sP}{{\mathcal P}}
\newcommand{\sQ}{{\mathcal Q}}
\newcommand{\sR}{{\mathcal R}}
\newcommand{\sS}{{\mathcal S}}
\newcommand{\sT}{{\mathcal T}}
\newcommand{\sU}{{\mathcal U}}
\newcommand{\sV}{{\mathcal V}}
\newcommand{\sW}{{\mathcal W}}
\newcommand{\sX}{{\mathcal X}}
\newcommand{\sY}{{\mathcal Y}}
\newcommand{\sZ}{{\mathcal Z}}
\newcommand{\ssA}{{\mathfrak A}}
\newcommand{\ssB}{{\mathfrak B}}
\newcommand{\ssC}{{\mathfrak C}}
\newcommand{\ssD}{{\mathfrak D}}
\newcommand{\ssE}{{\mathfrak E}}
\newcommand{\ssF}{{\mathfrak F}}
\newcommand{\ssg}{{\mathfrak g}}
\newcommand{\ssH}{{\mathfrak H}}
\newcommand{\ssI}{{\mathfrak I}}
\newcommand{\ssJ}{{\mathfrak J}}
\newcommand{\ssK}{{\mathfrak K}}
\newcommand{\ssL}{{\mathfrak L}}
\newcommand{\ssM}{{\mathfrak M}}
\newcommand{\ssN}{{\mathfrak N}}
\newcommand{\ssO}{{\mathfrak O}}
\newcommand{\ssP}{{\mathfrak P}}
\newcommand{\ssQ}{{\mathfrak Q}}
\newcommand{\ssR}{{\mathfrak R}}
\newcommand{\ssS}{{\mathfrak S}}
\newcommand{\ssT}{{\mathfrak T}}
\newcommand{\ssU}{{\mathfrak U}}
\newcommand{\ssV}{{\mathfrak V}}
\newcommand{\ssW}{{\mathfrak W}}
\newcommand{\ssX}{{\mathfrak X}}
\newcommand{\ssY}{{\mathfrak Y}}
\newcommand{\ssZ}{{\mathfrak Z}}

\newcommand{\B}{{\mathbb B}}
\newcommand{\C}{{\mathbb C}}
\newcommand{\D}{{\mathbb D}}
\newcommand{\E}{{\mathbb E}}
\newcommand{\F}{{\mathbb F}}
\newcommand{\G}{{\mathbb G}}
\renewcommand{\H}{{\mathbb H}}
\newcommand{\J}{{\mathbb J}}
\newcommand{\K}{{\mathbb K}}
\renewcommand{\L}{{\mathbb L}}
\newcommand{\M}{{\mathbb M}}
\newcommand{\N}{{\mathbb N}}
\renewcommand{\P}{{\mathbb P}}
\newcommand{\Q}{{\mathbb Q}}
\newcommand{\R}{{\mathbb R}}

\newcommand{\T}{{\mathbb T}}
\newcommand{\U}{{\mathbb U}}
\newcommand{\V}{{\mathbb V}}
\newcommand{\W}{{\mathbb W}}
\renewcommand{\S}{{\mathbb S}}
\newcommand{\Y}{{\mathbb Y}}
\newcommand{\Z}{{\mathbb Z}}
\newcommand{\id}{{  id}}

\newcommand{\A}{\"{a}}
\newcommand{\rank}{{rank}}
\newcommand{\END}{{\mathbb E}{ nd}}
\newcommand{\End}{{ End}}
\newcommand{\Hg}{{Hg}}
\newcommand{\tr}{{ tr}}
\newcommand{\Tr}{{ Tr}}
\newcommand{\SL}{{ SL}}
\newcommand{\PSL}{{ PSL}}
\newcommand{\Gl}{{ Gl}}
\newcommand{\Cor}{{ Cor}}
\newcommand{\Gal}{{ Gal}}
\newcommand{\GL}{{ GL}}
\newcommand{\PGL}{{ PGL}}
\newcommand{\MT}{{ MT}}
\newcommand{\Hdg}{{  Hdg}}
\newcommand{\MTV}{{  MTV}}
\newcommand{\SO}{{  SO}}
\newcommand{\Sp}{{ Sp}}
\newcommand{\Hom}{{ Hom}}
\newcommand{\Ker}{{ Ker}}
\newcommand{\Lie}{{ Lie}}
\newcommand{\Aut}{{ Aut}}
\newcommand{\Image}{{ Image}}
\newcommand{\Gr}{{ Gr}}
\newcommand{\Id}{{Id}}
\newcommand{\rk}{{ rk}}
\newcommand{\pardeg}{{par.deg}}
\newcommand{\SU}{{ SU}}
\newcommand{\Res}{{ Res}}
\newcommand{\Fr}{{ Frob_p}}
\newcommand{\Spec}{{ Spec}}
\newcommand{\Ext}{{ Ext}}
\newcommand{\Sym}{{ Sym}}
\newcommand{\Tor}{{Tor}}
\newcommand{\ch}{{ ch}}
\newcommand{\qtq}[1]{\quad\mbox{#1}\quad}
\newcommand{\bp}{\bar{\partial}}
\newcommand{\Om}{\Omega}
\newcommand{\td}{ Td}
\newcommand{\ind}{Ind}
\newcommand{\ds}{\oplus}
\newcommand{\bds}{\bigoplus}
\newcommand{\ts}{\otimes}
\newcommand{\bts}{\bigotimes}
\newcommand{\diag}{ diag}
\newcommand{\st}{\stackrel}
\newcommand{\btheorem}{\begin{theorem}}
\newcommand{\etheorem}{\end{theorem}}
\newcommand{\bproposition}{\begin{proposition}}
\newcommand{\eproposition}{\end{proposition}}
\newcommand{\bdefinition}{\begin{definition}}
\newcommand{\edefinition}{\end{definition}}
\newcommand{\bcorollary}{\begin{corollary}}
\newcommand{\ecorollary}{\end{corollary}}
\newcommand{\bproof}{\begin{proof}}
\newcommand{\eproof}{\end{proof}}
\newcommand{\bremark}{\begin{remark}}
\newcommand{\eremark}{\end{remark}}
\newcommand{\eexample}{\end{example}}
\newcommand{\bexample}{\begin{example}}
\newcommand{\la}{\langle}
\newcommand{\elemma}{\end{lemma}}
\newcommand{\blemma}{\begin{lemma}}
\newcommand{\ra}{\rangle}
\newcommand{\sq}{\sqrt{-1}}
\newcommand{\suml}{\sum\limits}
\newcommand{\jk}{dz_{j}\wedge d\bar{z}_{k}}
\newcommand{\ckm}{compact K\"{a}hler manifold\ }
\newcommand{\km}{K\"{a}hler manifold\ }
\newcommand{\hvb}{holomorphic vector bundle\ }
\newcommand{\hhvb}{hermitian holomorphic vector bundle\ }
\newcommand{\hlb}{holomorphic line bundle\ }
\newcommand{\hm}{hermitian manifold\ }
\newcommand{\hpd}{hermitian positive definite\ }
\newcommand{\KM}{K\"{a}hler metric\ }
\newcommand{\ob}{orthonormal basis\ }
\newcommand{\ccm}{compact complex manifold\ }
\newcommand{\hhlb}{hermitian holomorphic line bundle\ }
\newcommand{\chm}{compact hermitian manifold\ }
\newcommand{\ep}{\epsilon}
\newcommand{\om}{\omega}
\newcommand{\Llr}{\Longleftrightarrow}
\newcommand{\Lr}{\Longrightarrow}
\newcommand{\ldo}{linear differential operator\ }
\newcommand{\Supp}{ Supp}
\newcommand{\dsl}{\bigoplus\limits}
\newcommand{\wt}{\widetilde}

\newcommand{\Psh}{ Psh}
\newcommand{\sumo}{\widetilde{\suml}}
\newcommand{\grad}{ grad}
\newcommand{\emb}{\hookrightarrow}
\newcommand{\ut}{\la A_{E,\omega}u, u\ra}
\newcommand{\gt}{\la A^{-1}_{E,\omega}g, g\ra}
\newcommand{\pq}{\wedge^{p,q}T^{*}X\ts E}
\newcommand{\nq}{\wedge^{n,q}T^{*}X\ts E}
\newcommand{\pn}{\wedge^{p,n}T^{*}X\ts E}

\newcommand{\p}{\partial}
\newcommand{\pic}{ Pic}
\newcommand{\ord}{ord}
\newcommand{\Div}{ Div}

\renewcommand{\bar}{\overline}
\newcommand{\eps}{\varepsilon}
\newcommand{\pa}{\partial}
\renewcommand{\phi}{\varphi}
\newcommand{\OO}{\mathcal O}

\newcommand{\ka}{K\"ahler }
\newcommand{\kar}{K\"ahler-Ricci}
\newcommand{\ee}{\end{eqnarray*}}
\newcommand{\be}{\begin{eqnarray*}}
\newcommand{\red}{\textcolor[rgb]{1.00,0.00,0.00}}
\newcommand{\yel}{\textcolor[rgb]{0.00,0.00,1.00}}

\newcommand{\rt}{\right)}
\newcommand{\lt}{\left(}

\newcommand{\aver}[1]{-\hskip-0.35cm\int_{#1}}
\newcommand{\beq}{\begin{equation}}
\newcommand{\eeq}{\end{equation}}

\newcommand{\bd}{\begin{description}}
\newcommand{\ed}{\end{description}}
\newcommand{\ul}{\underline}
\renewcommand{\hat}{\widehat}
\renewcommand{\tilde}{\widetilde}

\newcommand{\rw}{\rightarrow}
\renewcommand{\bf}{\textbf}
\renewcommand{\sc}{\textsc}
\renewcommand{\it}{\textit}
\newcommand{\md}{\textmd}
\renewcommand{\sf}{\textsf}
\renewcommand{\rm}{\textrm}
\newcommand{\lrw}{\Longrightarrow}

\newcommand{\var}{\varnothing}

\newcommand{\dra}{\dashrightarrow}
\renewcommand{\bf}{\textbf}

\newcommand{\TO}{\Longrightarrow}
\newcommand{\OT}{\Longleftarrow}
\newcommand{\col}{\textcolor[rgb]{0.80,0.10,0.40}}
\renewcommand{\Hom}{\text{Hom}}
\newcommand{\wge}{{\text{\Large{$\wedge$}}}}
\newcommand{\pxj}{\frac{\p}{\p x^j}}
\newcommand{\pxi}{\frac{\p}{\p x^i}}
\newcommand{\pxk}{\frac{\p}{\p x^k}}
\newcommand{\pxl}{\frac{\p}{\p x^l}}

\newcommand{\pzj}{\frac{\p}{\p z^j}}
\newcommand{\pzi}{\frac{\p}{\p z^i}}
\newcommand{\pzk}{\frac{\p}{\p z^k}}
\newcommand{\pzl}{\frac{\p}{\p z^\ell}}
\newcommand{\pzs}{\frac{\p}{\p z^s}}
\newcommand{\pzt}{\frac{\p}{\p z^t}}
\newcommand{\pzp}{\frac{\p}{\p z^p}}
\newcommand{\pzq}{\frac{\p}{\p z^q}}
\newcommand{\pzA}{\frac{\p}{\p z^A}}
\newcommand{\pzB}{\frac{\p}{\p z^B}}
\newcommand{\pzC}{\frac{\p}{\p z^C}}
\newcommand{\pzD}{\frac{\p}{\p z^D}}

\newcommand{\bpzj}{\frac{\p}{\p \bar z^j}}
\newcommand{\bpzi}{\frac{\p}{\p \bar z^i}}
\newcommand{\bpzk}{\frac{\p}{\p \bar z^k}}
\newcommand{\bpzl}{\frac{\p}{\p \bar z^\ell}}
\newcommand{\bpzs}{\frac{\p}{\p \bar z^s}}
\newcommand{\bpzt}{\frac{\p}{\p \bar z^t}}
\newcommand{\bpzp}{\frac{\p}{\p \bar z^p}}
\newcommand{\bpzq}{\frac{\p}{\p \bar z^q}}

\newcommand{\ci}{I_i}
\newcommand{\cj}{I_j}
\newcommand{\ck}{I_k}
\newcommand{\cl}{I_\ell}
\newcommand{\cp}{I_p}
\newcommand{\cq}{I_q}
\newcommand{\cs}{I_s}
\newcommand{\ct}{I_t}

\newcommand{\cbi}{I_{ \bar i}}
\newcommand{\cbj}{I_{ \bar j}}
\newcommand{\cbk}{I_{\bar k}}
\newcommand{\cbl}{I_{\bar\ell}}
\newcommand{\cbp}{I_{ \bar p}}
\newcommand{\cbq}{I_{ \bar q}}
\newcommand{\cbs}{I_{\bar s}}
\newcommand{\cbt}{I_{ \bar t}}

\setlinespacing{1.0}
\def\baselinestretch{1.0}
{

\section{Introduction}

It is well-known(\cite{BS}) that on a compact K\"ahler manifold, if
the Ricci curvature is positive, then the first Betti number is
zero; if the Ricci curvature is negative, then there is no
holomorphic vector field. The key ingredient for the proofs of such
results is the K\"ahler symmetry. On the other hand, on an Hermitian
manifold, we don't have such symmetry and there are several
different Ricci curvatures. While on a K\"ahler manifold, all these
Ricci curvatures coincide, since the Chern curvature on a K\"ahler
manifold coincides with the curvature of the (complexified)
Levi-Civita connection. We can see this more clearly on
 an abstract Hermitian holomorphic bundle $E$. The Chern connection
 $\nabla^{CH}$  on $E$  is the unique connection which is compatible with
 the holomorphic structure and the Hermitian metric on $E$. Hence, the Chern
 curvature  $\Theta^{E}\in \Gamma(M, \Lambda^{1,1}T^*M\ts E^*\ts
 E)$. There are two methods to take trace of $\Theta^E$. If we take
 trace of $\Theta^E$ on the part $End(E)=E^*\ts E$, we get a $(1,1)$-form on
 $M$ which
 it is called the first Ricci curvature. It is well
 known that the first Ricci curvature represents the first Chern
 class of the bundle. On the other
 hand, if we take trace on the $(1,1)$-part using the metric of the
 manifold, we obtain  an endomorphism of $E$, $Tr_\omega \Theta^E \in \Gamma(M, E^*\ts
 E)$. It is called the second Ricci curvature of $\Theta^E$. The
 first and second Ricci curvatures
 have different geometric meanings, which were not clearly studied in some earlier literatures.
  We should point out that the nonexistence of holomorphic
 sections is characterized by the second Ricci curvature.
  Let $E$ be the holomorphic
 tangent bundle $T^{1,0}M$. If $M$ is K\"ahler,  the first and second Ricci curvature
 are the same by the K\"ahler symmetry. Unfortunately, on an
 Hermitian manifold, the Chern curvature is not symmetric, i.e.,
 the first and second Ricci curvatures are different. Moreover, in general they
 can not be compared. An interesting example is the Hopf
manifold $\S^{2n+1}\times \S^1$. The canonical metric on it has
strictly positive second Ricci curvature!

 In this paper, we study the nonexistence of holomorphic and harmonic sections of an
 abstract vector bundle over a compact Hermitian manifold. Let $E$
 be a holomorphic vector bundle over a compact Hermitian
 manifold $(M,\omega)$. Since the holomorphic section space $H^0(M,E)$ is
 independent on the connections of $E$, we can
 choose any connection on it. As we mentioned above, the key part, is the second
 Ricci curvature of the connection.  For example, on the holomorphic tangent bundle $T^{1,0}M$ of
  an Hermitian manifold $M$, there are three common connections
\bd
\item[(1)]  the complexified Levi-Civita connection $\nabla$ on
 $T^{1,0}M$;
\item[(2)]  the Chern connection $\nabla^{CH}$ on  $T^{1,0}M$;

\item[(3)] the Bismut connection $\nabla^B$
on  $T^{1,0}M$. \ed

\noindent It is well-known  that if $M$ is K\"ahler, all three
connections are the same. However,  in general, the relations among
them are somewhat mysterious. In this paper, we derive certain
relations about their curvatures on certain Hermitian manifolds.

Let $E$ be an  Hermitian \emph{complex} (possibly
\emph{non-holomorphic}) vector bundle or a Riemannian \emph{real}
vector bundle over a compact Hermitian manifold $(M,\omega)$. Let
$\p_E, \bp_E$ be the $(1,0), (0,1)$ part of $\nabla^E$ respectively.
The $(1,1)$-curvature of $\nabla^E$ is denoted by $R^E\in
\Gamma(M,\Lambda^{1,1}T^*M\ts E^*\ts E)$. It can be viewed as a
representation of the operator $\p_E\bp_E+\bp_E\p_E$.
 We can define harmonic
section spaces associated to $(E,\nabla^E)$ by \beq
\sH^{p,q}_{\bp_E}(M,E)=\{\phi\in \Om^{p,q}(M,E)\ |\
\bp_E\phi=\bp_E^*\phi=0\} \eeq In general, on a complex vector
bundle $E$, there is no such terminology like ``holomorphic section
of $E$". However, if the vector bundle $E$ is holomorphic
 and  $\nabla^E$ is the Chern connection on $E$ i.e. $\bp_E=\bp$, then $\sH^{p,q}_{\bp_E}(M,E)$ is
isomorphic to the Dolbeault cohomology group $H_{\bp}^{p,q}(M,E)$
and $H^{0}_{\bp}(M,E)$ is the holomorphic section space $H^{0}(M,E)$
of $E$.

\btheorem Let $E$ be an Hermitian complex vector bundle or a
Riemannian real vector bundle over a compact Hermitian manifold
$(M,\omega)$ and $\nabla^E$ be any metric connection on
$E$.\bd\item[(1)] If the second Hermitian-Ricci curvature $
Tr_\omega R^E$ is nonpositive everywhere, then every $\bp_E$-closed
section of $E$ is parallel, i.e. $\nabla^Es=0$;

\item[(2)] If  the second Hermitian-Ricci curvature $ Tr_\omega R^E$ is nonpositive everywhere and negative at some point, then $\sH^0_{\bp_E}(M,E)=0$;

\item[(3)]  If  the second Hermitian-Ricci curvature $ Tr_\omega R^E$ is $p$-nonpositive everywhere and $p$-negative at some point, then $\sH^0_{\bp_E}(M,\Lambda^qE)=0$ for any $p\leq q\leq
rank(E)$. \ed
 \etheorem

The proof of this theorem is based on generalized Bochner-Kodaira
identities on  vector bundles over Hermitian manifolds (Theorem
\ref{bochner}). We prove that (Theorem \ref{Gauduchonbochner}) the
torsion integral of the Hermitian manifold can be killed if the
background Hermitian metric $\omega$ is Gauduchon, i.e.
$\p\bp\omega^{n-1}=0$. On the other hand, in the conformal class of
any Hermitian metric, the Gauduchon metric always exists
(\cite{Ga4}). So we can change the background metric in the
conformal way and the positivity of the second Hermitian-Ricci
curvature is preserved. This method is very useful on Hermitian
manifolds. Kobayashi-Wu(\cite{KW}) and Gauduchon(\cite{Ga2})
obtained similar result in the special case when $\nabla^E$ is the
Chern connection of the Hermitian \emph{holomorphic} vector bundle
$E$. Now we go back to the Hermitian manifold $(M,\omega)$.

\bcorollary\label{A} Let $(M,\omega)$ be a compact Hermitian
manifold\bd
\item[(1)] if the second Ricci-Chern curvature $Tr_\omega\Theta$ is
nonnegative everywhere and positive at some point, then
$H^{p,0}_{\bp}(M)=0$ for any $1\leq p\leq n$. In particular, the
arithmetic genus $\chi(M,\sO)=1$;
\item[(2)] if
the second Ricci-Chern curvature $Tr_\omega\Theta$ is nonpositive
everywhere and negative at some point, then the holomorphic vector
bundle $\Lambda^{p}T^{1,0}M$ has no holomorphic vector field for any
$1\leq p\leq n$. \ed \ecorollary \noindent Since the first
Ricci-Chern curvature and the second Ricci-Chern curvature of an
Hermitian manifold can not be compared, we can not derive that the
manifold $M$ is K\"ahler, even if the second Ricci-Chern curvature
is positive everywhere. In general, the first Ricci-Chern curvature
is $d$-closed but the second Ricci-Chern curvature is not $d$-closed
and so they are in the different $(d,\bp,\p)$-cohomology classes.
For example, the Hopf manifold $\S^{2n+1}\times \S^1$ with standard
Hermitian metric has strictly positive second Ricci-Chern curvature
and nonnegative first Ricci-Chern curvature, but it is non-K\"ahler.
For more details, see Proposition \ref{hopf}.

%

Now we consider several special Hermitian manifolds. An interesting
class of  Hermitian manifolds is the balanced Hermitian manifolds,
i.e.,
 Hermitian manifolds with coclosed K\"ahler forms. It is well-known
 that
 every K\"ahler manifold is balanced. In some literatures, they are also called semi-K\"ahler manifolds.
  In complex dimension $1$ and  $2$,
 every balanced Hermitian manifold is K\"ahler. However, in
 higher dimensions, there exist non-K\"ahler manifolds which admit
 balanced Hermitian metrics. Such examples were constructed by
 E. Calabi(\cite{CE}), see also \cite{GA} and \cite{Mic}.
There are also some other
 important classes of non-K\"ahler balanced manifolds, such as: complex
solvmanifolds, 1-dimensional families of K\"ahler manifolds (see
\cite{Mic}) and compact complex parallelizable manifolds (except
complex torus) (see \cite{UH}). On the other hand, Alessandrini-
Bassaneli( \cite{AB}) proved that every Moishezon manifold is
balanced and so  balanced manifolds can be constructed from K\"ahler
manifolds by modification. For more examples, see \cite{AB3},
\cite{Mic}, \cite{GI1} and \cite{GI2}.

Every balanced metric $\omega$ is a Gauduchon metric. In fact,
$d^*\omega=0$ is equivalent to $d\omega^{n-1}=0$ and so
$\p\bp\omega^{n-1}=0$. By \cite{Ga4}, every Hermitian manifold has a
Gauduchon metric. However, there are many manifolds which can not
support balanced metrics.  For example, the Hopf surface $\S^3\times
\S^1$ is non-K\"ahler, so it has no balanced metric. For more
discussion , one can see \cite{CE}, \cite{Mic},\cite{TW2}, \cite{AB}
and \cite{AB3}.

On a compact balanced Hermitian manifold $M$, we can detect the
holomorphic section spaces $H^{p,0}_{\bp}(M)$ by   Levi-Civita
connection. Let $\nabla$ be the complexified Levi-Civita connection
and $\nabla'$, $\nabla''$ the $(1,0)$ and $(0,1)$ components of
$\nabla$ respectively. In general, holomorphic $p$-forms are not
$\nabla''$-closed. The Ricci curvatures related to the Levi-Civita
connection are defined in \ref{riccilc} and \ref{2riccilc}.

\btheorem\label{B} Let $(M,\omega)$ be a compact balanced Hermitian
manifold.
 If the Hermitian-Ricci curvature $(R_{i\bar j})$ of $M$ is nonnegative
everywhere, then \bd
\item[(1)] If $\phi$ is a holomorphic $p$-form, then
$\Delta_{\p}\phi=0$ and so $h^{p,0}(M)\leq h^{0,p}(M)$ for any
$1\leq p\leq n$;

\item[(2)] If the Hermitian-Ricci curvature $(R_{i\bar j})$ is
 positive at some point, then $
H^{p,0}_{\bp}(M)=0  $ for any $1\leq p\leq n$. In particular, the
arithmetic genus $\chi(M,\sO)=1$. \ed

 \etheorem

The dual of Theorem \ref{B} is

\btheorem\label{C} Let $(M,\omega)$ be a compact balanced Hermitian
manifold. If $2\hat R^{(2)}_{i\bar j}-R_{i\bar j}$ is nonpositive
everywhere and negative at some point, there is no holomorphic
vector field on M. \etheorem

\bremark  It is easy to see that the Hermitian-Ricci curvature
tensor $(R_{i\bar j})$ and second Ricci-Chern curvature tensor
$\Theta^{(2)}:=Tr_\omega\Theta$ can not be compared. Therefore,
Theorem \ref{B} and Corollary \ref{A} are independent of each other.
For the same reason, Theorem \ref{C} and Corollary \ref{A} are
independent. Balanced Hermitian manifolds with nonnegative
Hermitian-Ricci curvatures are discussed in Proposition
\ref{positivebalanced}. \eremark

As we discuss in the above, on  Hermitian manifolds, the second
Ricci curvature tensors of various metric connections are closely
related to the geometry of Hermitian manifolds. A natural idea is to
define a flow by using second Ricci curvature tensors of various
metric connections. For example, \beq \frac{\p h}{\p
t}=-\Theta^{(2)}+\mu h,\ \ \ \mu\in\R \eeq on a general Hermitian
manifold $(M,h)$ by using the second Ricci-Chern curvature. This
flow preserves the K\"ahler and the Hermitian structure and has
short time solution on any compact Hermitian manifold. It is very
similar to and closely related to the Hermitian Yang-Mills flow, the
K\"ahler-Ricci flow and the harmonic map heat flow. It may be a
bridge to connect them. In this paper we only briefly discuss its
basic properties. In a subsequent paper(\cite{LLY}) we will study
its geometric and analytic property in detail.

We would like to thank Yi Li, Jeffrey Streets, Valetino Tosatti for
their useful comments on an earlier version of this paper.

\section{Various connections and curvatures on Hermitian manifolds}

\subsection{Complexified Riemannian curvature}
Let $(M,g)$ be a Riemannian manifold with Levi-Civita connection
$\nabla$, the curvature $R$ of $(M,g,\nabla)$ is defined  as \beq
R(X,Y,Z,W)=g\left(\left(\nabla_X\nabla_Y-\nabla_Y\nabla_X-\nabla_{[X,Y]}\right)Z,W
\right)\eeq On an Hermitian manifold $(M,h)$, let $\nabla$ be the
complexified Levi-Civita connection and $g$ the background
Riemannian metric. Two metrics are related by  \beq
ds^2_h=ds^2_g-\sq \omega_h \eeq where $\omega_h$ is the fundamental
$(1,1)$-form (or K\"ahler form) associated to $h$. For any two
holomorphic vector fields $X,Y\in\Gamma(M, T^{1,0}M)$, \beq
h(X,Y)=2g(X,\bar Y) \eeq \emph{This formula will be used  in several
definitions}.
 In the local holomorphic coordinates
$\{z^1,\cdots, z^n\}$ on $M$, the complexified Christoffel symbols
are given by \beq
\Gamma_{AB}^C=\sum_{E}\frac{1}{2}g^{CE}\big(\frac{\p g_{AE}}{\p
z^B}+\frac{\p g_{BE}}{\p z^A}-\frac{\p g_{AB}}{\p
z^E}\big)=\sum_{E}\frac{1}{2}h^{CE}\big(\frac{\p h_{AE}}{\p
z^B}+\frac{\p h_{BE}}{\p z^A}-\frac{\p h_{AB}}{\p z^E}\big)
 \eeq
where $A,B,C,E\in \{1,\cdots,n,\bar{1},\cdots,\bar{n}\}$ and
$z^{A}=z^{i}$ if $A=i$, $z^{A}=\bar{z}^{i}$ if $A=\bar{i}$. For
example \beq \Gamma_{ij}^k=\frac{1}{2}h^{k\bar \ell}\left(\frac{\p
h_{j\bar \ell}}{\p z^i}+\frac{\p h_{i\bar \ell}}{\p z^j}\right),\
\Gamma_{\bar ij}^k=\frac{1}{2}h^{k\bar \ell}\left(\frac{\p h_{j\bar
\ell}}{\p \bar z^i}-\frac{\p h_{j\bar i}}{\p \bar z^\ell}\right)
\eeq The complexified curvature components are
 \be
 R_{ABCD}:&=&2\bf{g}\left(\left(\nabla_{\pzA}\nabla_{\pzB}-\nabla_{\pzB}\nabla_{\pzA}\right)\pzC,\pzD
\right)\\&=&
\bf{h}\left(\left(\nabla_{\pzA}\nabla_{\pzB}-\nabla_{\pzB}\nabla_{\pzA}\right)\pzC,\frac{\p}{\p
z^{\bar D}} \right)\ee Hence \beq
R_{ABC}^D=\sum_ER_{ABCE}h^{ED}=-\left(\frac{\p\Gamma_{AC}^D}{\p
z^B}-\frac{\p\Gamma_{BC}^D}{\p
z^A}+\Gamma_{AC}^F\Gamma_{FB}^D-\Gamma_{BC}^F\Gamma_{AF}^D\right)
\eeq  By the Hermitian property, we have, for example \beq
 R_{i\bar jk}^{l}=-\left(\frac{\p
\Gamma^{l}_{ik}}{\p \bar z^j}-\frac{\p \Gamma^{l}_{\bar jk}}{\p
z^i}+\Gamma_{ ik}^{s}\Gamma^{l}_{\bar js}-\Gamma_{ \bar
jk}^{s}\Gamma^{l}_{ is}-{\Gamma_{\bar j k}^{\bar s}\Gamma_{i\bar
s}^l}\right)\eeq
%
%

%
%
%

\bremark We have $R_{ABCD}=R_{CDAB}$. In particular, \beq R_{i\bar
jk\bar\ell}=R_{k\bar \ell i\bar j} \eeq \eremark

\noindent Unlike the K\"ahler case, we can define several Ricci
curvatures:

\bdefinition \bd\item[(1)] The \emph{complexified Ricci curvature}
on $(M, h)$ is defined by \beq \mathscr{R}_{k\bar \ell}:=h^{i\bar
j}\left(R_{k\bar j i\bar \ell}+R_{ki\bar j \bar\ell}\right)
\label{cr}\eeq The \emph{complexified scalar curvature} of $h$ is
defined as \beq s_h:=h^{k\bar\ell}\mathscr{R}_{k\bar \ell}\eeq
\item[(2)] The \emph{Hermitian-Ricci curvature}  is  \beq R_{k\bar
\ell}:=h^{i\bar j}R_{i\bar j k\bar \ell}\label{riccilc} \eeq
 The \emph{Hermitian-scalar curvature} of $h$ is given by\beq
S:=h^{k\bar \ell}R_{k\bar \ell} \eeq \ed\edefinition

\blemma On an Hermitian manifold, \beq \bar{R_{ABCD}}=R_{\bar A\bar
B\bar C\bar D},\ \ \bar {\mathscr{R}_{k\bar
\ell}}=\mathscr{R}_{\ell\bar k},\ \ \ \bar{R_{k\bar
\ell}}=R_{\ell\bar k} \eeq and \beq \mathscr R_{k\bar \ell}=h^{i\bar
j}\left(2R_{k\bar j i \bar\ell}-R_{k\bar \ell i\bar j}\right) \eeq
\bproof The Hermitian property of curvature tensors is obvious. By
first Bianchi identity, we have
$$ R_{k
i\bar j \bar\ell}+R_{k\bar j \bar\ell i}+R_{k\bar \ell i\bar j}=0 $$

\noindent That is $ R_{ki\bar j \bar\ell}=R_{k\bar j
i\bar\ell}-R_{k\bar \ell i\bar j} $. The curvature formula \ref{cr}
turns to be \beq \mathscr{R}_{k\bar \ell}=h^{i\bar j}\left(2R_{k\bar
j i\bar \ell}-R_{k\bar \ell i\bar j}\right) \eeq

\eproof
 \elemma

\bdefinition The Ricci curvatures are called \emph{positive} ( resp.
\emph{nonnegative, negative, non-positive}) if the corresponding
Hermitian matrices are positive ( resp. nonnegative, negative,
non-positive). \edefinition

\noindent The following three formulas are used frequently in the
sequel. \blemma Assume $h_{i\bar j}=\delta_{ij}$ at a fixed point
$p\in M$, we have the following formula
\begin{eqnarray} R_{i\bar j k\bar
\ell}\nonumber&=&-\frac{1}{2}\left(\frac{\p^2 h_{i\bar \ell}}{\p
z^k\p\bar z^j}+\frac{\p^2 h_{k\bar j}}{\p z^i\p\bar
z^\ell}\right)\\
\nonumber &&+\frac{1}{4}\left(\frac{\p h_{k\bar q}}{\p z^i}\frac{\p
h_{q\bar\ell}}{\p\bar z^j}+\frac{\p h_{i\bar q}}{\p z^k}\frac{\p
h_{q\bar j}}{\p\bar z^\ell}\right)+\frac{1}{4}\left(\frac{\p
h_{i\bar q}}{\p z^k}\frac{\p h_{q\bar\ell}}{\p\bar z^j}+\frac{\p
h_{k\bar q}}{\p z^i}\frac{\p h_{q\bar j}}{\p\bar
z^\ell}\right)\\
&&+\frac{1}{4}\left(\frac{\p h_{q\bar \ell}}{\p z^i}\frac{\p
h_{k\bar j}}{\p\bar z^q}+\frac{\p h_{q\bar j}}{\p z^k}\frac{\p
h_{i\bar \ell}}{\p\bar z^q}\right) +
    \frac{1}{4}\left(\frac{\p h_{i\bar \ell}}{\p
z^q}\frac{\p h_{k\bar q}}{\p\bar z^j}+\frac{\p h_{k\bar j}}{\p
z^q}\frac{\p h_{i\bar q}}{\p\bar z^\ell}\right)\\
\nonumber&& -\frac{1}{4}\left(\frac{\p h_{q\bar \ell}}{\p
z^i}\frac{\p h_{k\bar q}}{\p\bar z^j}+\frac{\p h_{q\bar j}}{\p
z^k}\frac{\p h_{i\bar q}}{\p\bar
z^\ell}\right)-\frac{1}{4}\left(\frac{\p h_{i\bar \ell}}{\p
z^q}\frac{\p h_{k\bar j}}{\p\bar z^q}+\frac{\p h_{k\bar j}}{\p
z^q}\frac{\p h_{i\bar \ell}}{\p\bar z^q}\right)
\end{eqnarray}

\elemma

By a linear transformation on the local holomorphic coordinates, one
can get the following Lemma. For more details, we refer the reader
to \cite{ST2}. \blemma\label{normalcoordinate} Let $(M,h,\omega)$ be
an Hermitian manifold. For any $p\in M$, there exist local
holomorphic coordinates $\{z^i\}$ centered at a point $p$ such that
\beq h_{i\bar j}(p)=\delta_{ij} \qtq{and} \Gamma_{ij}^k(p)=0 \eeq
\elemma

 \noindent By  Lemma \ref{normalcoordinate},  we have a simplified version of curvatures:
\blemma  Assume $h_{i\bar j}(p)=\delta_{ij}$ and
$\Gamma_{ij}^k(p)=0$ at a fixed point $p\in M$, \beq R_{i\bar j
k\bar \ell}=-\frac{1}{2}\left(\frac{\p^2 h_{i\bar \ell}}{\p
z^k\p\bar z^j}+\frac{\p^2 h_{k\bar j}}{\p z^i\p\bar
z^\ell}\right)-\sum_q\left(\frac{\p h_{q\bar \ell}}{\p z^i}\frac{\p
h_{k\bar q}}{\p\bar z^j}+\frac{\p h_{q\bar j}}{\p z^k}\frac{\p
h_{i\bar q}}{\p\bar z^\ell}\right) \eeq For Hermitian-Ricci
curvatures \beq R_{k\bar \ell}=h^{i\bar j} R_{i\bar jk\bar
\ell}=-\frac{1}{2}\sum_s\left(\frac{\p^2 h_{s\bar \ell}}{\p
z^k\p\bar z^s}+\frac{\p^2 h_{k\bar s}}{\p z^s\p\bar
z^\ell}\right)-\sum_{q,s}\left(\frac{\p h_{q\bar \ell}}{\p
z^s}\frac{\p h_{k\bar q}}{\p\bar z^s}+\frac{\p h_{k\bar q}}{\p
z^s}\frac{\p h_{q\bar \ell}}{\p\bar
z^s}\right)\label{generalricci1}\eeq and \beq h^{i\bar j} R_{ k\bar
j i\bar \ell}=h^{i\bar j}R_{i\bar \ell k\bar
j}=-\frac{1}{2}\sum_s\left(\frac{\p^2 h_{k\bar \ell}}{\p z^s\p\bar
z^s}+\frac{\p^2 h_{s\bar s}}{\p z^k\p\bar
z^\ell}\right)-\sum_{q,s}\left(\frac{\p h_{q\bar \ell}}{\p
z^k}\frac{\p h_{s\bar q}}{\p \bar z^s}+\frac{\p h_{q\bar s}}{\p
z^s}\frac{\p h_{k\bar q}}{\p \bar z^\ell}\right) \eeq For
complexified Ricci curvature,
\begin{eqnarray} \mathscr R_{k\bar
\ell}&=&\nonumber\frac{1}{2}\sum_s\left(\frac{\p^2 h_{s\bar
\ell}}{\p z^k\p\bar z^s}+\frac{\p^2 h_{k\bar s}}{\p z^s\p\bar
z^\ell}\right)- \sum_s\left(\frac{\p^2 h_{k\bar \ell}}{\p z^s\p\bar
z^s}+\frac{\p^2
h_{s\bar s}}{\p z^k\p\bar z^\ell}\right)\\
&+& \sum_{q,s}\left(\frac{\p h_{q\bar \ell}}{\p z^s}\frac{\p
h_{k\bar q}}{\p\bar z^s}+\frac{\p h_{k\bar q}}{\p z^s}\frac{\p
h_{q\bar \ell}}{\p\bar z^s}\right)-2\sum_{q,s}\left(\frac{\p
h_{q\bar \ell}}{\p z^k}\frac{\p h_{s\bar q}}{\p \bar z^s}+\frac{\p
h_{q\bar s}}{\p z^s}\frac{\p h_{k\bar q}}{\p \bar z^\ell}\right)
\label{curvatureformula}\end{eqnarray}  \elemma

\subsection{Curvature of complexified Levi-Civita connection on $T^{1,0}M$}

Since $T^{1,0}M$ is a subbundle of $T_{\C}M$, there is an induced
connection $\hat\nabla$ on $T^{1,0}M$ given by \beq
\hat\nabla=\pi\circ\nabla:
T^{1,0}M\stackrel{\nabla}{\rightarrow}\Gamma(M, T_{\C}M\ts
T_{\C}M)\stackrel{\pi}{\rightarrow}\Gamma(M,T_{\C}M\ts T^{1,0}M)
\label{metricconnection}\eeq The curvature $\hat R\in
\Gamma(M,\Lambda^2 T_{\C}M\ts T^{*1,0}M\ts T^{1,0}M)$ of
$\hat\nabla$ is given by \beq \hat R(X,Y)s
=\hat\nabla_{X}\hat\nabla_Ys-\hat\nabla_Y\hat\nabla_Xs-\hat\nabla_{[X,Y]}s\eeq
for any $X,Y\in T_{\C}M$ and $s\in T^{1,0}M$. It has components \beq
\hat R_{ABk}^l=\frac{\p\Gamma_{Bk}^l}{\p
z^A}-\frac{\p\Gamma_{Ak}^{l}}{\p
z^B}-\Gamma_{Ak}^{s}\Gamma_{Bs}^{l}+\Gamma_{Bk}^{s}\Gamma_{As}^{l}\eeq
where \beq \hat R\left(\frac{\p}{\p z^A}, \frac{\p}{\p
z^B}\right)\frac{\p}{\p z^k}=\sum_{l}\hat R_{ABk}^l\pzl\eeq
 For example,
\beq \hat R_{i\bar jk}^{l}=-\left(\frac{\p \Gamma^{l}_{ik}}{\p \bar
z^j}-\frac{\p \Gamma^{l}_{\bar jk}}{\p z^i}+\Gamma_{
ik}^{s}\Gamma^{l}_{\bar js}-\Gamma_{ \bar jk}^{s}\Gamma^{l}_{s
i}\right)\label{11curvature}\eeq
%

\noindent With respect to the Hermitian metric $h$ on $T^{1,0}M$, we
can define \beq \hat R_{AB k\bar l}=\sum_{s=1}^n\hat R_{AB
k}^{s}h_{s\bar \ell} \eeq

\bdefinition The \emph{first Ricci curvature} of the Hermitian
vector bundle $\left(T^{1,0}M,\hat \nabla\right)$ is defined by \beq
\hat R^{(1)}_{i\bar j}=h^{k\bar \ell} \hat R_{i\bar j k\bar\ell}
\eeq The \emph{second Ricci curvature} of it is
 \beq \hat R^{(2)}_{k\bar \ell}=h^{i\bar j} \hat R_{i\bar j k\bar
 \ell}\label{2riccilc}
\eeq The \emph{scalar curvature} of $\hat\nabla$ on $T^{1,0}M$ is
denoted by \beq S^{LC}=h^{i\bar j}h^{k\bar \ell}\hat R_{i\bar j
k\bar \ell} \eeq

\edefinition

\noindent By Lemma \ref{normalcoordinate}, we have the following
formulas

\blemma On an Hermitian manifold $(M,h)$, on a point $p$ with
$h_{i\bar j}(p)=\delta_{ij}$ and $\Gamma_{ij}^k(p)=0$,
\begin{eqnarray} \hat R_{i\bar j k\bar
\ell}=-\frac{1}{2}\left(\frac{\p^2 h_{i\bar \ell}}{\p z^k\p\bar
z^j}+\frac{\p^2 h_{k\bar j}}{\p z^i\p\bar
z^\ell}\right)-\sum_q\frac{\p h_{q\bar \ell}}{\p z^i}\frac{\p
h_{k\bar q}}{\p\bar z^j}
\end{eqnarray}
For the Ricci curvatures,
 \beq \hat R^{(1)}_{i\bar
j}= -\frac{1}{2}\sum_k\left(\frac{\p^2 h_{i\bar k}}{\p z^k\p\bar
z^j}+\frac{\p^2 h_{k\bar j}}{\p z^i\p\bar
z^k}\right)-\sum_{k,q}\frac{\p h_{q\bar k}}{\p z^i}\frac{\p h_{k\bar
q}}{\p\bar z^j} \label{1ricci}\eeq and \beq \hat R^{(2)}_{i\bar
j}=-\frac{1}{2}\sum_k\left(\frac{\p^2 h_{i\bar k}}{\p z^k\p\bar
z^j}+\frac{\p^2 h_{k\bar j}}{\p z^i\p\bar
z^k}\right)-\sum_{k,q}\frac{\p h_{i\bar q}}{\p\bar z^k}\frac{\p
h_{q\bar j}}{\p z^k}\label{2ricci}\eeq
Moreover, \beq \hat R^{(1)}_{i\bar j}-\hat R_{i\bar
j}^{(2)}=h_{m\bar j}h^{\ell\bar k}\Gamma_{\bar k i}^{\bar
q}\Gamma_{\ell \bar q}^m-\Gamma_{k\bar j}^{\bar q}\Gamma_{i\bar
q}^k=\sum_{k,q}\left(\frac{\p h_{i\bar q}}{\p\bar z^k}\frac{\p
h_{q\bar j}}{\p z^k}-\frac{\p h_{i\bar q}}{\p z^k}\frac{\p h_{q\bar
j}}{\p \bar z^k} \right)\eeq

 \elemma

%

\subsection{Curvature of Chern connection on $T^{1,0}M$}
On the Hermitian holomorphic vector bundle $(T^{1,0}M,h)$, the Chern
connection $\nabla^{CH}$ is the unique connection which is
compatible with the complex structure and the Hermitian metric. Its
curvature components are \beq \Theta_{i\bar j k\bar
\ell}=-\frac{\p^2 h_{k\bar \ell}}{\p z^i\p \bar z^j}+h^{p\bar
q}\frac{\p h_{p\bar \ell}}{\p \bar z^j}\frac{\p h_{k\bar q}}{\p
z^i}\label{cherncurvature} \eeq
 It is
well-known that the \emph{first Ricci-Chern curvature} \beq
\Theta^{(1)}:= \frac{\sq}{2\pi}\Theta^{(1)}_{i\bar j} dz^i\wedge
d\bar z^j\eeq represents the first Chern class of $M$ where
  \beq
\Theta^{(1)}_{i\bar j}= h^{k\bar \ell}\Theta_{i\bar j k\bar \ell}
=-\frac{\p^2 \log \det(h_{k\bar \ell})}{\p z^i\p\bar
z^j}\label{1chern}\eeq The \emph{second Ricci-Chern curvature}
components are\beq \Theta^{(2)}_{i\bar j}=h^{k\bar
\ell}\Theta_{k\bar \ell i\bar j} \label{2chern}\eeq The \emph{scalar
curvature} of the Chern connection is defined by \beq
S^{CH}=h^{i\bar j}h^{k\bar \ell}\Theta_{i\bar j k\bar \ell}\eeq

\subsection{Curvature of Bismut connection on  $T^{1,0}M$}

 In \cite{bismut}, Bismut defined a class of connections on Hermitian
 manifolds. In this
 subsection, we choose one of them (see \cite{MM}, p. $21$).
  The \emph{Bismut connection} $\nabla^B$ on the holomorphic tangent bundle
 $(T^{1,0}M,h)$ is characterized by
 \beq \nabla^B=\nabla+S^B \eeq
where $S^B$ is  $1$-form with values in $End(T^{1,0}M)$\beq
\bf{h}(S^B(X)Y, Z)= 2\bf{g}(S^B(X)Y,\bar Z)=\sq
(\p-\bp)\omega_h(X,Y,\bar Z) \eeq for any $Y,Z\in T^{1,0} M$ and
$X\in T_{\C} M$. Let $\tilde \Gamma_{i\alpha}^\beta$ and $\tilde
\Gamma_{\bar j\alpha}^\beta$ be the Christoffel symbols of the
Bismut connection where $i,j,\alpha,\beta\in\{1,\cdots, n\}$. We use
different types of letters since the Bismut connection is not
torsion free.

 \blemma\label{bismutlc} We have the
following relations between $\tilde \Gamma$ and $\Gamma$, \beq
\tilde\Gamma_{i\alpha\bar\beta}: =h_{\beta\bar\gamma
}\Gamma_{i\alpha}^{\bar
\gamma}=\Gamma_{i\alpha\bar\beta}+\Gamma_{\alpha\bar\beta
i}=\frac{\p h_{i\bar\beta}}{\p z^\alpha},\ \ \ \ \tilde\Gamma_{\bar
j\alpha\bar\beta}=2\Gamma_{\bar
j\alpha\bar\beta}\label{symbolrelations}\eeq \bproof Let $X=\pzi,
Y=\pzj, Z=\pzk$. Since $\omega_h=\frac{\sq}{2}h_{m\bar n}dz^m\wedge
d\bar z^n$, we obtain \be \sq (\p-\bp)\omega_h(X,Y,\bar
Z)&=&-\frac{1}{2}\frac{\p h_{m\bar n}}{\p z^p}dz^pdz^md\bar z^n
\left(\pzi,\pzj,\bpzk\right)\\&=&\frac{1}{2}\left(\frac{\p h_{i\bar
k}}{\p z^j}-\frac{\p h_{j\bar k}}{\p z^i}\right)\\&=&\Gamma_{j\bar
k}^{\bar s}h_{i\bar s}=\Gamma_{j\bar k i} \ee On the other hand \beq
h\left(\nabla^B_{\pzi}\pzj,\pzk\right)=\tilde \Gamma_{ij\bar k} \eeq
Using the definition of Bismut connection, we get \beq
\tilde\Gamma_{i\alpha\bar\beta}
=\Gamma_{i\alpha\bar\beta}+\Gamma_{\alpha\bar\beta i}=\frac{\p
h_{i\bar\beta}}{\p z^\alpha} \eeq The proof of the other one is
similar. \eproof \elemma

  The Bismut curvature $B\in\Gamma\left(M,\Lambda^{1,1}T^*M\ts
End(T^{1,0}M)\right)$ is given by \beq B_{i\bar
j\alpha}^\beta=-\frac{\p \tilde\Gamma_{i\alpha}^\beta}{\p\bar
z^j}+\frac{\p\tilde \Gamma_{\bar j\alpha}^\beta}{\p
z^i}-\tilde\Gamma_{i\alpha}^\gamma\tilde\Gamma_{\bar
j\gamma}^\beta+\tilde\Gamma_{\bar
j\alpha}^\gamma\tilde\Gamma_{i\gamma}^\beta
\label{curvaturedefinition}\eeq

\blemma Assume $h_{i\bar j}(p)=\delta_{ij}$ and $\Gamma_{ij}^k(p)=0$
at a fixed point $p\in M$, \beq B_{i\bar j
\alpha\bar\beta}=-\left(\frac{\p^2 h_{i\bar\beta}}{\p\bar z^j\p
z^\alpha}+\frac{\p^2 h_{\alpha\bar j}}{\p z^i\p\bar
z^\beta}-\frac{\p^2 h_{\alpha\bar\beta}}{\p z^i\p\bar
z^j}\right)+\sum_\gamma\frac{\p h_{\alpha\bar\gamma}}{\p
z^i}\frac{\p h_{\gamma\bar\beta}}{\p\bar z^j}-4\sum_\gamma\frac{\p
h_{\alpha\bar\gamma}}{\p \bar z^j}\frac{\p h_{\gamma\bar\beta}}{\p
z^i}\label{bismutcurvature} \eeq \bproof It follows by
\ref{symbolrelations} and \ref{curvaturedefinition}.\eproof
 \elemma

\noindent We can define  the first Ricci-Bismut curvature
$B^{(1)}_{i\bar j}$, the second Ricci-Bismut curvature
$B^{(2)}_{i\bar j}$ and scalar curvature $S^{BM}$ similarly.

\subsection{Relations among the four curvatures on Hermitian manifolds}

\bproposition\label{curvaturecompare} On an Hermitian manifold
$(M,h)$, we have \beq R_{ijk\bar l}=\hat R_{ijk\bar \ell},\ \ \
R_{\bar i\bar j k\bar \ell}=\hat R_{\bar i\bar j k\bar \ell} \eeq
 and for any $u,v\in \C^n$, \beq\left( R_{i\bar j k\bar \ell}-\hat  R_{i\bar j k\bar
\ell}\right) u^{i}\bar u^j v^k\bar v^\ell\leq 0\eeq In particular,
$R_{i\bar j}\leq \hat R^{(1)}_{i\bar j}$ and $R_{i\bar j}\leq \hat
R^{(2)}_{i\bar j}$ in the sense of Hermitian matrices. \bproof Let
\beq T_{i\bar j k\bar \ell}= R_{i\bar j k\bar \ell}-\hat R_{i\bar j
k\bar \ell}=\Gamma_{\bar j k}^{\bar s}\Gamma_{i\bar s}^{ t} h_{t\bar
\ell} \eeq Without loss generality, we assume $h_{i\bar
j}=\delta_{ij}$ at a fixed point, then \beq T_{i\bar j k\bar
\ell}=\sum_s\Gamma_{\bar j k s}\Gamma_{i\bar s \bar \ell}=-\sum_s
\Gamma_{i\bar s\bar \ell}\bar{\Gamma_{j\bar s \bar k}}\eeq where
\beq \Gamma_{i\bar s \bar \ell}=\frac{1}{2}\left(\frac{\p h_{i\bar
\ell}}{\p\bar z^s}-\frac{\p h_{i\bar s}}{\p \bar
z^\ell}\right)=-\Gamma_{i\bar\ell\bar s} \eeq and so $T_{i\bar j
k\bar \ell}u^{i}\bar u^j v^k\bar v^\ell\leq 0$. \eproof
\eproposition

\bremark\label{cohomology} \bd\item[(1)] Because of the second order
terms in $R$, $\hat R$, $\Theta$ and $B$, we can not compare $R,\hat
R$ with $\Theta$, $B$.
\item[(2)] Since the third order terms of $\p \Theta^{(2)}$ are not
zero in general. Therefore it is possible that $\Theta^{(1)}$ and
$\Theta^{(2)}$ are not in the same  $(d,\p,\bp)$-cohomology class.
For the same reason $B^{(1)}$ and $B^{(2)}$ are not in the same
$(d,\p,\bp)$-cohomology class.
\item[(3)] If the manifold $(M,h)$ is K\"ahler, then all curvatures
are the same. \ed\eremark

\section{Curvature relations on special Hermitian manifolds}

\subsection{Curvatures relations on balanced Hermitian manifolds}

\noindent The following lemma is well-known( for example
\cite{Ga1}), and we include a proof here in our setting.
 \blemma\label{balancedef} Let $(M,\omega)$ be a compact Hermitian manifold. The following
 conditions are equivalent:

 \bd\item[(1)] $d^*\omega=0$;

 \item[(2)] $d\omega^{n-1}=0$;

 \item[(3)] For any smooth function $f\in C^\infty(M)$,
\beq \frac{1}{2}\Delta_d f=\Delta_{\bp}f=\Delta_\p f=-h^{i\bar
j}\frac{\p^2 f}{\p z^i\p\bar z^j} \eeq

\item[(4)] $ \Gamma_{\bar i\ell}^\ell=0 $ for any $1\leq i\leq n$.
\ed \bproof On a compact Hermitian manifold,
$d^*\omega=-*d*\omega=-c_n*d\omega^{n-1}$ where $c_n$ is a constant
depending only on the complex dimension $n$ of $M$. On the other
hand, the Hodge $*$ is an isomorphism, and so $(1)$ and $(2)$ are
equivalent. If $f$ is a smooth function on $M$, \beq
\begin{cases} \Delta_{\bp} f=-h^{i\bar j}\frac{\p^2 f}{\p z^i\p\bar
z^j}+2h^{i\bar j}\Gamma_{i\bar j}^{\bar \ell}\frac{\p
f}{\p\bar z^\ell}\\
\Delta_\p f=-h^{i\bar j}\frac{\p^2 f}{\p z^i\p\bar z^j}+2h^{i\bar
j}\Gamma_{\bar j i}^{k}\frac{\p f}{\p z^k} \end{cases}\eeq On the
other hand, \beq h^{i\bar j}\Gamma_{i\bar j}^{\bar
\ell}=-\Gamma_{k\bar j}^{\bar j}h^{k\bar \ell} \qtq{and} h^{i\bar
j}\Gamma_{\bar j i}^{k}=-\Gamma_{\bar \ell i}^ih^{k\bar \ell}\eeq
Therefore $(3)$ and $(4)$ are equivalent. For the equivalence of
$(1)$ and $(4)$, see Lemma \ref{balancedcondition}.
 \eproof
 \elemma

\bdefinition An Hermitian manifold $(M,\omega)$ is called
\emph{balanced} if  it satisfies one of the conditions in Lemma
\ref{balancedef}. \edefinition

\noindent On a balanced Hermitian manifold, there are more
symmetries on the second derivatives of the metric.\blemma Let
$(M,h)$ be a balanced Hermitian manifold. On a point $p$ with
$h_{i\bar j}(p)=\delta_{ij}$ and $\Gamma_{ij}^k(p)=0$, \beq
\sum_{s}\frac{\p h_{s\bar i }}{\p\bar  z^s} = \sum_s \frac{\p
h_{s\bar s}}{\p \bar z^i}=0\label{mm}\eeq and \beq\sum_i \frac{\p^2
h_{i\bar \ell} }{\p z^k\p\bar z^i}=\sum_i\frac{\p^2 h_{k\bar i} }{\p
z^i\p\bar z^\ell}=\sum_i\frac{\p^2 h_{i\bar i}}{\p z^k\p\bar
z^\ell}-2\sum_{i,q}\frac{\p h_{q\bar \ell}}{\p\bar z^i}\frac{\p
h_{k\bar q}}{\p z^i} \label{pp}\eeq \bproof At a fixed point $p$, if
$h_{i\bar j}=0$ and $\Gamma_{ij}^k=0$, then \beq \frac{\p h_{i\bar
j}}{\p\bar z^k}=-\frac{\p h_{i\bar k}}{\p\bar z^j}\label{123} \eeq
The balanced condition $\sum_s\Gamma_{\bar is}^s=0$ is reduced to
\beq \sum_s\frac{\p h_{s\bar s}}{\p \bar z^i}=\sum_s \frac{\p
h_{s\bar i}}{\p\bar z^s}=0 \eeq by formula \ref{123}. By the
balanced condition \be 0=\frac{\p \Gamma_{\bar \ell i}^i}{\p
z^k}&=&\frac{\p}{\p z^k}\left(\frac{1}{2} h^{i\bar q}\left(\frac{\p
h_{i\bar q}}{\p \bar z^\ell}-\frac{\p h_{i\bar \ell}}{\p\bar
z^q}\right)\right)\\&=& \frac{1}{2}\sum_i\left( \frac{\p^2 h_{i\bar
i}}{\p z^k\p\bar z^\ell}-\frac{\p^2 h_{i\bar \ell} }{\p z^k\p\bar
z^i}\right)- \sum_{i,q}\frac{\p h_{q\bar \ell}}{\p\bar z^i}\frac{\p
h_{k\bar q}}{\p z^i}\ee Hence, we obtain formula \ref{pp}.
 \eproof

\elemma

\bproposition\label{balancericci} Let $(M,h)$ be a balanced
Hermitian manifold. At a point $p$ with $h_{i\bar j}(p)=\delta_{ij}$
and $\Gamma_{ij}^k(p)=0$, we have following formulas about various
Ricci curvatures:
\begin{eqnarray} \Theta^{(1)}_{k\bar \ell}&=&\hat R^{(1)}_{k\bar
\ell}=B^{(1)}_{k\bar \ell}=-\sum_i\frac{\p^2 h_{i\bar i}}{\p z^k
\p\bar z^\ell}+\sum_{q,i}\frac{\p h_{q\bar \ell}}{\p\bar
z^i}\frac{\p h_{k\bar q}}{\p z^i}\label{2} \\ \Theta^{(2)}_{k\bar
\ell}&=&-\sum_i\frac{\p^2 h_{k\bar \ell}}{\p z^i \p\bar
z^i}+\sum_{i,q}\frac{\p h_{q\bar \ell}}{\p\bar z^i}\frac{\p h_{k\bar
q}}{\p z^i}\label{complexifiedricci2}\\ \hat R^{(2)}_{k\bar
\ell}&=&-\sum_i\frac{\p^2 h_{i\bar i}}{\p z^k \p\bar
z^\ell}+\sum_{i,q}\left(2\frac{\p h_{q\bar \ell}}{\p\bar
z^i}\frac{\p h_{k\bar q}}{\p z^i}-\frac{\p h_{k\bar q}}{\p\bar
z^i}\frac{\p h_{q\bar \ell}}{\p z^i}\right) \\
B_{k\bar \ell}^{(2)}&=&-\sum_i\frac{\p^2 h_{i\bar i}}{\p z^k \p\bar
z^\ell}+\sum_{i,q}\left(5\frac{\p h_{q\bar \ell}}{\p\bar
z^i}\frac{\p h_{k\bar q}}{\p z^i}-4\frac{\p h_{k\bar q}}{\p\bar
z^i}\frac{\p h_{q\bar \ell}}{\p z^i}\right)
\\
 R_{k\bar
\ell}&=&-\sum_i\frac{\p^2 h_{i\bar i}}{\p z^k \p\bar
z^\ell}+\sum_{i,q}\left(\frac{\p h_{q\bar \ell}}{\p\bar z^i}\frac{\p
h_{k\bar q}}{\p z^i}-\frac{\p h_{k\bar q}}{\p\bar z^i}\frac{\p
h_{q\bar \ell}}{\p z^i}\right)\label{3} \\ \mathscr{R}_{k\bar
\ell}&=&- \sum_i \frac{\p^2 h_{k\bar \ell}}{\p z^i\p\bar
z^i}-\sum_{i,q}\left(\frac{\p h_{q\bar \ell}}{\p\bar z^i}\frac{\p
h_{k\bar q}}{\p z^i}-\frac{\p h_{k\bar q}}{\p\bar z^i}\frac{\p
h_{q\bar \ell}}{\p z^i}\right)\label{ricciflowcurvature}
\end{eqnarray}

\bproof In  \ref{1ricci}, \ref{2ricci}, \ref{1chern},\ref{2chern},
\ref{generalricci1}, \ref{curvatureformula}, we get expressions for
all Ricci curvatures on Hermitian manifolds. By balanced relations
\ref{mm} and \ref{pp}, we get simplified versions of all Ricci
curvatures. \eproof
 \eproposition

 \bproposition\label{positivebalanced} \bd\item[(1)] A balanced Hermitian manifold
with positive Hermitian-Ricci curvature $R_{i\bar j}$ is K\"ahler.

\item[(2)] Let $(M,h)$ be a compact balanced Hermitian manifold. If the
 Hermitian-Ricci curvature is nonnegative
everywhere and positive at some point, then $M$
  is Moishezon. \ed \bproof
\bf{(1)} On a balanced Hermitian manifold \beq \Theta^{(1)}_{i\bar
j}=\hat R^{(1)}_{i\bar j}\geq R_{i\bar j}\eeq If $R_{i\bar j}$ is
Hermitian positive, then $\Theta^{(1)}_{i\bar j}$ is Hermitian
positive, and so \beq \Om =-\frac{\sq}{2\pi}\p\bp\log\det(h_{k\bar
\ell})\eeq is a K\"ahler metric.

\bf{(2)} If the
 Hermitian-Ricci curvature is nonnegative
everywhere and positive at some point, so is $\Theta^{(1)}_{i\bar
j}$. The Hermitian line bundle $L=\det(T^{1,0}M)$ satisfies \beq
\int_M c_1(L)^n> 0\eeq By Siu-Demailly's solution of
Grauert-Riemenschneider conjecture (\cite{Siu} \cite{D4}), $M$ is
Moishezon. \eproof
 \eproposition

\subsection{Curvature relations on Hermitian manifolds with $\Lambda(\p\bp\omega)=0$}

Now we consider a compact Hermitian manifold $(M,\omega)$ with
$\Lambda(\p\bp\omega)=0$. The condition $\Lambda(\p\bp\omega)=0$ is
equivalent to \beq \sum_k\left(\frac{\p h_{i\bar j}}{\p z^k \p\bar
z^k}+\frac{\p h_{k\bar k}}{\p z^i\p\bar
z^j}\right)=\sum_k\left(\frac{\p h_{i\bar k}}{\p z^k\p\bar
z^j}+\frac{\p^2 h_{k\bar j}}{\p z^i\p\bar z^k}\right) \label{yy}\eeq
for any $i,j$. We can use  \ref{yy} to simplify Ricci curvatures and
get relations among them.

\bproposition\label{ddbarricci} Let $(M,h)$ be a compact Hermitian
manifold with $\Lambda(\p\bp\omega)=0$. At a point $p$ with
$h_{i\bar j}(p)=\delta_{ij}$ and $\Gamma_{ij}^k(p)=0$, the following
identities about Ricci curvatures hold:
\begin{eqnarray} \Theta^{(1)}_{k\bar \ell}&=&-\sum_i\frac{\p^2 h_{i\bar i}}{\p z^k
\p\bar z^\ell}+\sum_{q,i}\frac{\p h_{q\bar \ell}}{\p\bar
z^i}\frac{\p h_{k\bar q}}{\p z^i}\label{2} \\ \Theta^{(2)}_{k\bar
\ell}&=&-\sum_i\frac{\p^2 h_{k\bar \ell}}{\p z^i \p\bar
z^i}+\sum_{i,q}\frac{\p h_{q\bar \ell}}{\p\bar z^i}\frac{\p h_{k\bar
q}}{\p z^i}\label{complexifiedricci21}\\
\hat R^{(1)}_{k\bar \ell}&=&-\frac{1}{2}\sum_i\left(\frac{\p^2
h_{k\bar\ell}}{\p z^i\p\bar z^i }+\frac{\p^2 h_{i\bar i}}{\p z^k
\p\bar z^\ell}\right)-\sum_{i,q}\frac{\p h_{q\bar \ell}}{\p\bar
z^i}\frac{\p h_{k\bar q}}{\p z^i}
 \\\hat R^{(2)}_{k\bar
\ell}&=&-\frac{1}{2}\sum_i\left(\frac{\p^2 h_{k\bar\ell}}{\p
z^i\p\bar z^i }+\frac{\p^2 h_{i\bar i}}{\p z^k \p\bar
z^\ell}\right)-\sum_{i,q}\frac{\p h_{k\bar q}}{\p\bar
z^i}\frac{\p h_{q\bar \ell}}{\p z^i} \\
B^{(1)}_{k\bar\ell}&=&-\sum_i\frac{\p^2 h_{k\bar \ell}}{\p z^i\p\bar
z^i}+\sum_{i,q}\left(\frac{\p h_{q\bar \ell}}{\p\bar z^i}\frac{\p
h_{k\bar q}}{\p z^i}-4\frac{\p h_{q\bar \ell}}{\p z^i}\frac{\p
h_{k\bar q}}{\p\bar z^i}\right)\\
B_{k\bar \ell}^{(2)}&=&-\sum_i\frac{\p^2 h_{i\bar i}}{\p z^k \p\bar
z^\ell}+\sum_{i,q}\left(\frac{\p h_{q\bar \ell}}{\p\bar z^i}\frac{\p
h_{k\bar q}}{\p z^i}-4\frac{\p h_{k\bar q}}{\p\bar z^i}\frac{\p
h_{q\bar \ell}}{\p z^i}\right)
\\
 R_{k\bar
\ell}&=&-\frac{1}{2}\sum_i\left(\frac{\p^2 h_{k\bar\ell}}{\p
z^i\p\bar z^i }+\frac{\p^2 h_{i\bar i}}{\p z^k \p\bar
z^\ell}\right)-\sum_{i,q}\left(\frac{\p h_{q\bar \ell}}{\p\bar
z^i}\frac{\p h_{k\bar q}}{\p z^i}+\frac{\p h_{k\bar q}}{\p\bar
z^i}\frac{\p h_{q\bar \ell}}{\p z^i}\right)\label{3} \\
\mathscr{R}_{k\bar \ell}&=&-\frac{1}{2}\sum_i\left(\frac{\p^2
h_{k\bar\ell}}{\p z^i\p\bar z^i }+\frac{\p^2 h_{i\bar i}}{\p z^k
\p\bar z^\ell}\right)+\sum_{i,q}\left(\frac{\p h_{q\bar
\ell}}{\p\bar z^i}\frac{\p h_{k\bar q}}{\p z^i}+\frac{\p h_{k\bar
q}}{\p\bar z^i}\frac{\p h_{q\bar \ell}}{\p
z^i}\right)\\\nonumber&&-2\sum_{q,i}\left(\frac{\p h_{q\bar
\ell}}{\p z^k}\frac{\p h_{i\bar q}}{\p \bar z^i}+\frac{\p h_{q\bar
i}}{\p z^i}\frac{\p h_{k\bar q}}{\p \bar z^\ell}\right)
\end{eqnarray}
\eproposition

\bproposition\label{chernbismut} If $(M,\omega)$ is a compact
Hermitian manifold with $\Lambda(\p\bp\omega)=0$, then \beq
B^{(2)}\leq \Theta^{(1)} \qtq{and} B^{(1)}\leq \Theta^{(2)} \eeq in
the sense of Hermitian matrices and identities hold if and only if
$(M,\omega)$ is K\"ahler. Moreover, \beq
\Theta^{(2)}+B^{(2)}=\Theta^{(1)}+R^{(1)} \eeq

  \eproposition

Finally, we would like to discuss the relations of balanced
manifolds and strong K\"ahler manifolds with torsion. By \cite{AB},
every Moishezon manifold is balanced, i.e. there exists a smooth
Hermitian metric $\omega$ such that $d^*\omega=0$. On the other
hand, by Demailly-Paun \cite{DP}( see also \cite{JS}),  on each
Moishezon manifold, there exists a singular Hermitian metric
$\omega$ such that $\p\bp\omega=0$ in the sense of current. However,
these two conditions can not be satisfied simultaneously in the
smooth sense on an Hermitian non-K\"ahler manifold. It is known in
\cite{AI} and also \cite{FPS}, but merits a proof in our setting.

 \bproposition Let $(M,\omega)$ be a compact Hermitian manifold. If
$d^*\omega=0$ and $\Lambda(\p\bp\omega)=0$, then $d\omega=0$, i.e.
$(M,\omega)$ is K\"ahler. In particular, if a compact Hermitian
manifold admits a smooth metric $\omega$ such that $d^*\omega=0$ and
$\p\bp\omega=0$, then it is K\"ahler. \bproof Let $(M,\omega)$ be a
balanced Hermitian manifold with $\Lambda(\p\bp\omega)=0$. The
condition $\Lambda(\p\bp\omega)=0$ is equivalent to \beq
\sum_i\frac{\p h_{i\bar i}}{\p z^k\p\bar z^\ell}+\sum_i\frac{\p^2
h_{k\bar \ell}}{\p z^i\p\bar z^i}= \sum_i\frac{\p h_{i\bar \ell}}{\p
z^k\p\bar z^i}+\sum_i\frac{\p^2 h_{k\bar i}}{\p z^i \p\bar
z^\ell}\eeq By formula \ref{pp}, at a point $p$ with $h_{i\bar
j}=\delta_{ij}$ and $\Gamma_{ij}^k(p)=0$, we have \be \sum_i\frac{\p
h_{i\bar i}}{\p z^k\p\bar z^\ell}+\sum_i\frac{\p^2 h_{k\bar
\ell}}{\p z^i\p\bar z^i}&=&\sum_i \frac{\p h_{i\bar \ell}}{\p
z^k\p\bar
z^i}+\sum_i\frac{\p^2 h_{k\bar i}}{\p z^i \p\bar z^\ell}\\
&=&2\sum_i\frac{\p h_{i\bar i}}{\p z^k\p\bar
z^\ell}-4\sum_{q,i}\frac{\p h_{q\bar \ell}}{\p\bar z^i}\frac{\p
h_{k\bar q}}{\p z^i}\ee That is \beq \sum_i\frac{\p h_{i\bar i}}{\p
z^k\p\bar z^\ell}=\sum_i\frac{\p^2 h_{k\bar \ell}}{\p z^i\p\bar
z^i}+4\sum_{q,i}\frac{\p h_{q\bar \ell}}{\p\bar z^i}\frac{\p
h_{k\bar q}}{\p z^i}\eeq Taking trace of it, we obtain \beq
4\sum_{q,i, k}\frac{\p h_{q\bar k}}{\p\bar z^i}\frac{\p h_{k\bar
q}}{\p z^i}=0 \Longleftrightarrow\frac{\p h_{k\bar q}}{\p z^i}=0\eeq
at point $p$. Since $p$ is arbitrary, we have $d\omega\equiv0$,
therefore, $(M,\omega)$ is K\"ahler. \eproof
 \eproposition

\section{Bochner formulas on Hermitian complex and Riemannian real vector bundles over compact Hermitian manifolds}

Let $(M,h,\omega)$ be a compact Hermitian manifold. The complexified
Levi-Civita connection $\nabla$  on $T_{\C}M$ induces a linear
connection on $\Om^{p,q}(M)$: \beq \nabla:
\Om^{p,q}(M)\>>>\Om^1(M)\ts \left(\Om^{p,q}(M)\ds\Om^{p-1,q+1}(M)\ds
\Om^{p+1,q-1}(M)\right) \eeq We consider the following two canonical
components of $\nabla$, \beq\begin{cases}
\nabla':\Om^{p,q}(M)\rightarrow \Om^{1,0}(M)\ts \Om^{p,q}(M)\\
\nabla'':\Om^{p,q}(M)\rightarrow \Om^{0,1}(M)\ts
\Om^{p,q}(M)\end{cases}\eeq Note that $\nabla\neq \nabla'+\nabla''$
if $(M,h,\omega)$ is not K\"ahler. The following
 calculation rule follows immediately \beq\nabla'(\phi\wedge
\psi)=\left(\nabla'\phi\right)\wedge \psi+\phi\wedge \nabla'\psi
\eeq for any $\phi,\psi\in \Om^{\bullet}(M)$.

\blemma On an Hermitian manifold $(M,h)$, we have \beq\begin{cases}
\p
h(\phi,\psi)=h(\nabla'\phi,\psi)+h(\phi,\nabla''\psi)\\
 \bp
h(\phi,\psi)=h(\nabla''\phi,\psi)+h(\phi,\nabla'\psi)
\end{cases} \qquad \Longleftrightarrow \begin{cases} \pzi
h(\phi,\psi)=h(\nabla'_i\phi,\psi)+h(\phi,\nabla''_{\bar i}\psi)\\
 \bpzj
h(\phi,\psi)=h(\nabla''_{\bar j}\phi,\psi)+h(\phi,\nabla'_{j}\psi)
\end{cases}
\eeq for any $\phi,\psi\in \Om^{p,q}(M)$. \elemma

\bremark \bd\item[(1)]
 Here we use the compact
notations
$$\nabla'_i=\nabla'_{\frac{\p}{\p z^i}},\ \ \nabla''_{\bar j}=\nabla''_{\frac{\p}{\p \bar z^j}}$$
Note that $\nabla'_{\bar j}=\nabla''_i=0$ and $\nabla_i\neq
\nabla'_i$, $\nabla_{\bar j}\neq \nabla'_{\bar j}$.
\item[(2)] If we regard $\Lambda^{p,q}T^*M$ as an abstract vector
bundle $E$, the above lemma says that $\nabla'$ and $\nabla''$ are
compatible with the Hermitian metric on $E$. \ed\eremark

Now we go to an abstract setting. Let $E$ be an  Hermitian
\emph{complex} (possibly \emph{non-holomorphic}) vector bundle or a
\emph{Riemannian} real vector bundle over a compact Hermitian
manifold $(M,\omega)$. There is a natural decomposition \beq
\nabla=\nabla^{'E}+\nabla{''^E} \eeq where \beq\begin{cases}
\nabla^{'E}:\Gamma(M,E)\>>>\Om^{1,0}(M,E)\\
\nabla^{''E}:\Gamma(M,E)\>>>\Om^{0,1}(M,E)
  \end{cases}\eeq
$\nabla^{'E}$ and $\nabla^{''E}$  induce two differential operators.
The first one is
 $
\p_E:\Om^{p,q}(M,E)\>>>\Om^{p+1,q}(M,E) $ defined by \beq
\p_E(\phi\ts s)=\left(\p\phi\right)\ts s+(-1)^{p+q}\phi\wedge
\nabla^{'E}s \eeq for any $\phi\in \Om^{p,q}(M)$ and $s\in
\Gamma(M,E)$. The other one is $
\bp_E:\Om^{p,q}(M,E)\>>>\Om^{p+1,q}(M,E) $ defined by \beq
\bp_E(\phi\ts s)=\left(\bp\phi\right)\ts s+(-1)^{p+q}\phi\wedge
\nabla^{''E}s \eeq for any $\phi\in \Om^{p,q}(M)$ and $s\in
\Gamma(M,E)$. The following formula is well-known \beq
\left(\p_E\bp_E+\bp_E\p_E\right)(\phi\ts s)=\phi\wedge
\left(\p_E\bp_E+\bp_E\p_E\right)s\eeq for any $\phi\in \Om^{p,q}(M)$
and $s\in \Gamma(M,E)$. The operator $\p_E\bp_E+\bp_E\p_E$ is
represented by its $(1,1)$ curvature tensor  $R^E \in \Gamma(M,
\Lambda^{1,1}T^*M\ts E)$. For any $\phi, \psi\in
\Om^{\bullet,\bullet}(M,E)$, there is a \emph{sesquilinear pairing}
\beq \left\{ \phi, \psi\right\}=\phi^\alpha\wedge \bar {\psi^\beta}
\langle e_\alpha, e_\beta\rangle \eeq  if $\phi=\phi^\alpha
e_\alpha$ and $\psi=\psi^\beta e_\beta$ in the local frames
$\{e_\alpha\}$ on $E$. By the metric compatible property of
$\nabla^E$, \beq \p\{\phi,\psi\}=\{\p_E
\phi,\psi\}+(-1)^{p+q}\{\phi, \bp_E\psi\} \eeq if
$\phi\in\Om^{p,q}(M,E)$.

 Let $\omega$ be the K\"ahler
form of the Hermitian metric $h$, i.e., \beq
\omega=\frac{\sq}{2}h_{i\bar j}dz^i\wedge d\bar z^j \eeq
 On the Hermitian manifold $(M,h,\omega)$, the norm on
 $\Om^{p,q}(M)$ is defined by
\beq
(\phi,\psi)=\int_M\la\phi,\psi\ra\frac{\omega^n}{n!}=\frac{2^{n}}{(p+q)!}
\int_M h(\phi,\psi)\frac{\omega^n}{n!}=\int_M\phi\wedge
*\bar \psi\label{norm}\eeq
 The norm on
$\Om^{p,q}(M,E)$ is defined by \beq (\phi,\psi)=\int_M
\{\phi,*\psi\}=\int_M\left(\phi^\alpha\wedge *\bar
{\psi^\beta}\right) \langle e_\alpha, e_\beta\rangle  \eeq for
$\phi,\psi\in \Om^{p,q}(M,E)$. The dual operators of $\p,\bp,\p_E$
and $\bp_E^*$ are denoted by $\p^*,\bp^*, \p^*_E$ and $\bp_E^*$
respectively.

The following lemma was firstly shown by Demailly using Taylor
expansion method( e.g. \cite{D}). For the convenience of the reader,
we will take another approach  which seems to be useful in local
computations.
 \blemma\label{Demailly} Let $(M,h,\omega)$  be a compact Hermitian manifold. If $\tau$ is the operator of type $(1,0)$ defined
by $\tau=[\Lambda,2\partial\omega]$ on $\Om^{\bullet}(M,E)$, \beq
\begin{cases}
\left[\Lambda,\p\right]=\sq \left(\bp^*+\bar\tau^*\right)\\
\left[\Lambda, \bp\right]=-\sq (\p^*+\tau^*)
\end{cases}
\eeq For the dual equation, it is \beq
\begin{cases}[\bp^*,L]=\sq (\p+\tau)\\
[\p^*,L]=-\sq(\bp+\bar\tau)
\end{cases}
\eeq  where $L$ is the operator $L\phi=2\omega\wedge\phi$ and
$\Lambda$ is
 the adjoint operator of $L$.

 \bproof See Appendix Lemma \ref{DemaillyA}. \eproof \elemma

In the rest of this section $E$ is assumed to be an Hermitian
complex vector bundles or a Riemannian real vector bundle over
 a compact Hermitian manifold $M$. \blemma\label{LY} Let $\nabla^E$ be a
metric connection on  $E$ over a compact Hermitian manifold
$(M,\omega)$. If $\tau$ is the operator of type $(1,0)$ defined by
$\tau=[\Lambda,2\partial\omega]$ on $\Om^{\bullet}(M,E)$, then\bd
\item[(1)] $[\bp_E^*,L]=\sq(\p_E+\tau)$;
\item[(2)] $[\p^*_E,L]=-\sq(\bp_E+\bar{\tau})$; \item[(3)]
$[\Lambda,\p_E]=\sqrt{-1}(\bp_E^*+\bar{\tau}^{*})$ ;
\item[(4)]
$[\Lambda,\bp_E]=-\sqrt{-1}(\p_E^*+\tau^{*})$.\ed

\bproof See Appendix Lemma \ref{LYA}.  \eproof \elemma

\btheorem\label{bochner}  Let $\nabla^E$ be a metric connection  $E$
over a compact Hermitian manifold $(M,\omega)$.

 \beq
\Delta_{\bp_E}=\Delta_{\p_E}+\sq\left[\p_E\bp_E+\bp_E\p_E,\Lambda\right]+(\p_E
\tau^*+\tau^*\p_E)-(\bp_E\bar\tau^*+\bar\tau^*\bp_E)\label{Laplaciancompare}\eeq
where \beq \begin{cases}\Delta_{\bp_E}=\bp_E\bp_E^*+\bp_E^*\bp_E\\
\Delta_{\p_E}=\p_E\p_E^*+\p_E^*\p_E\end{cases} \eeq\bproof It
follows from Lemma \ref{LY}. \eproof \etheorem

We make a useful observation on the torsion $\tau$:

 \blemma\label{tau1} For any
$s\in \Gamma(M,E)$, we have \beq \tau (s) =-2\sq
\left(\bp^*\omega\right) \cdot s,\ \ \ \ \ \  \bar\tau (s) =2\sq
\left(\p^*\omega\right) \cdot s \label{tau}\eeq

 \bproof By definition \be \left([\Lambda
,2\p\omega]\right)s&=&2 \Lambda\left( (\p\omega )\cdot
s\right)\\&=&2\left(\Lambda(\p\omega)\right)\cdot s\\&=&-2\sq
\left(\bp^*\omega\right)\cdot  s \ee Here we use the identity \beq
\bp^*\omega=\sq \Lambda(\p\omega) \eeq where the proof of it is
contained in Lemma \ref{balancedcondition} of the Appendix. \eproof
\elemma

\bcorollary\label{balancedlaplaciancompare} If $(M,\omega)$ is a
compact balanced Hermitian manifold, and $\nabla^E$  a metric
connection on  $E$ over $M$, then \beq \|\bp_E s\|^2=\|\p_E
s\|^2+\left(\sq\left[\p_E\bp_E+\bp_E\p_E,\Lambda\right]s,s\right)\eeq
for any $s\in\Gamma(M,E)$. \bproof Since for any $s\in \Gamma(M,E)$,
$\tau s=\bar \tau s=0$ and $\tau^*s=\bar\tau^* s=0$ on a balanced
Hermitian manifold. \eproof
 \ecorollary

\btheorem\label{Gauduchonbochner} Let $(M,\omega)$ be an Hermitian
manifold with $\p\bp\omega^{n-1}=0$. If $\nabla^E$ is a metric
connection on  $E$ over $M$, then \beq0= \|\bp_E s\|^2=\|\p_E
s\|^2+\left(\sq\left[\p_E\bp_E+\bp_E\p_E,\Lambda\right]s,s\right)\label{gauduchonlaplaciancompare}
\eeq for any $s\in \Gamma(M,E)$ with $\bp_Es=0$. \bproof We only
have to prove that \beq\left((\p_E
\tau^*+\tau^*\p_E)s-(\bp_E\bar\tau^*+\bar\tau^*\bp_E)s,s\right)=0
\eeq which is equivalent to $ \left(\p_Es,\tau s\right)=0 $ since
$\tau^*s=\bar\tau^*s=\bp_Es=0$. By formula \ref{tau} and Stokes'
Theorem, \be \left(\tau^*\p_Es, s\right)&=&\left(\p_E s, \tau
s\right)=\int_M\left\{ \p_E s, *(\tau s)\right\}\\&=&2\sq\int_M
\left\{ \p_E s, *\left(\bp^*\omega\cdot s
\right)\right\}\\
&=&2\sq \int_M\left\{ \p_E s,\left(*\bp^*\omega\right) \cdot s\right\}\\
&=& -2\sq \int_M\left\{  s,\bp_E\left(\left(*\bp^*\omega\right) \cdot s\right)\right\}\\
 &=&-2\sq \int_M\left\{ s, \left(\bp *\bp^*\omega\right)\cdot
s- \left(*\bp^*\omega\right)\wedge \bp_E s\right\}
 \ee
It is easy to see that \beq
\bp*\bp^*\omega=-\bp**\p*\omega=c_n\bp\p\omega^{n-1}=0 \eeq since
$*\omega=c_n\omega^{n-1}$ where $c_n$ is a constant depending only
on the complex dimension of $M$. Hence \beq (\p_Es,\tau s)=2\sq
\int_M \left\{s, \left(*\bp^*\omega\right)\wedge \bp_E s\right\}=0
\eeq since $\bp_E s=0$. \eproof \etheorem


\bremark By these formula, we can obtain classical vanishing
theorems on K\"ahler manifolds and  rigidity of harmonic maps
between Hermitian and Riemannian manifolds. \eremark

\section{Vanishing theorems on  Hermitian manifolds}

\subsection{Vanishing theorems on compact Hermitian manifolds}

Let $E$ be an  Hermitian \emph{complex} (possibly
\emph{non-holomorphic}) vector bundle or a {Riemannian} \emph{real}
vector bundle over a compact Hermitian manifold $(M,\omega)$. Let
$\p_E, \bp_E$ be the $(1,0), (0,1)$ part of $\nabla^E$ respectively.
The $(1,1)$-curvature of $\nabla^E$ is denoted by $R^E\in
\Gamma(M,\Lambda^{1,1}T^*M\ts E^*\ts E)$. It is a representation of
the operator $\p_E\bp_E+\bp_E\p_E$.
 We can define harmonic
section spaces associated to $(E,\nabla^E)$ by \beq
\sH^{p,q}_{\bp_E}(M,E)=\{\phi\in \Om^{p,q}(M,E)\ |\
\bp_E\phi=\bp_E^*\phi=0\} \eeq In general, on a complex vector
bundle $E$, there is no such terminology like ``holomorphic section
of $E$". However, if the vector bundle $E$ is holomorphic
 and  $\nabla^E$ is the Chern connection of $E$ i.e. $\bp_E=\bp$, then $\sH^{p,q}_{\bp_E}(M,E)$ is
isomorphic to the Dolbeault cohomology group $H_{\bp}^{p,q}(M,E)$
and $H^{0}_{\bp}(M,E)$ is the holomorphic section spaces
$H^{0}(M,E)$ of $E$.

\bdefinition Let $A$ be an $r\times r$ Hermitian matrix and
$\lambda_1\leq \cdots\leq \lambda_r$ be eigenvalues of $A$. $A$ is
said to be \emph{$p$-nonnegative} (resp. \emph{positive, negative,
nonpositive}) for $1\leq p\leq r$ if \beq \lambda_{i_1}+\cdots
+\lambda_{i_p}\geq 0 (\qtq{resp.}  >0, <0, \leq 0)\qtq{for any}
1\leq i_1<i_2<\cdots <i_p\leq n\eeq \edefinition

\btheorem \label{vanishing} Let $\nabla^E$ be any metric connection
of an Hermitian {complex} vector bundle or a Riemannian real vector
bundle $E$ over a compact Hermitian manifold $(M,h,\omega)$.
 \bd\item[(1)] If the second Hermitian-Ricci curvature $ Tr_\omega R^E$ is nonpositive everywhere, then every $\bp_E$-closed section of $E$ is parallel, i.e. $\nabla^Es=0$;

\item[(2)] If  the second Hermitian-Ricci curvature $ Tr_\omega R^E$ is nonpositive everywhere and negative at some point, then $\sH^0_{\bp_E}(M,E)=0$;

\item[(3)]  If  the second Hermitian-Ricci curvature $ Tr_\omega R^E$ is $p$-nonpositive everywhere and $p$-negative at some point, then $\sH^0_{\bp_E}(M,\Lambda^qE)=0$ for any $p\leq q\leq
rank(E)$. \ed

\bproof By \cite{Ga4}, there exists a smooth function $u:M\>>>\R$
such that $\omega_G=e^u \omega$ is a Gauduchon metric, i.e.
$\p\bp\omega^{n-1}_G=0$.  Now we replace the metric $\omega$ on $M$
by the Gauduchon metric $\omega_G$.  By the relation
$\omega_G=e^u\omega$, we get
 \beq Tr_{\omega_G}R^E=e^{-u}Tr_\omega R^E\eeq Therefore, the positivity conditions in
the Theorem are preserved. Let $s\in\Gamma(M,E)$ with $\bp_E s=0$,
by formula \ref{gauduchonlaplaciancompare}, we obtain
 \beq 0=\|\p_E s\|^2+\left(\sq\left[\p_E\bp_E+\bp_E\p_E,\Lambda_
G\right]s,s\right)= \|\p_E s\|^2-\left(Tr_{\omega_G} R^E
s,s\right)\label{zz}\eeq where
 \beq R^E=\p_E\bp_E+\bp_E\p_E=R_{i\bar
j\alpha}^{ \beta}dz^i\wedge d\bar z^j\ts e^\alpha\ts e_\beta \eeq
Since the second Hermitian-Ricci curvature $Tr_{\omega_G} R^E$ has
components \beq R_{\alpha\bar\beta}=h^{i\bar j}_GR_{i\bar
j\alpha\bar\beta} \eeq  formula \ref{zz} can be written as \beq
0=\|\p_E s\|^2-\int_MR_{\alpha\bar \beta} s^\alpha\bar
s^{\beta}\label{a2} \eeq Now \bf{(1)} and \bf{(2)} follow by
identity \ref{a2} with the curvature conditions immediately. For
$\bf{(3)}$, we set $F=\Lambda^q E$ with $p\leq q\leq
 r=
 rank(E)$. Let $\lambda_1\leq \cdots \leq \lambda_r$ be the eigenvalues
 of $Tr_{\omega_G} R^E$, then we know
 \beq \lambda_1+\cdots+\lambda_p\geq 0 \eeq
and it is strictly positive at some point. If $p\leq q\leq r$, the
smallest eigenvalue of $Tr_{\omega_G}R^F$ is $\lambda_1+\cdots
+\lambda_q\geq 0$ and it is strictly positive at some point. By
\bf{(2)}, we know $\sH^0_{\bp_E}(M,F)=0$.
 \eproof
\etheorem

If $\nabla^E$ is the Chern connection of the Hermitian holomorphic
vector bundle $E$, we know $$ \sH_{\bp_E}^0(M,E)\cong H^0(M,E)
$$ since $\bp_E=\nabla^{''E}=\bp$ for the Chern connection.

\bcorollary[Kobayashi-Wu\cite{KW}, Gauduchon \cite{Ga2}] Let
$\nabla^E$ be the Chern connection of an Hermitian holomorphic
vector bundle $E$ over a compact Hermitian manifold $(M,h,\omega)$.
 \bd\item[(1)] If the second Ricci-Chern curvature $ Tr_\omega R^E$ is nonpositive everywhere, then every holomorphic section of $E$ is parallel, i.e. $\nabla^Es=0$;

\item[(2)] If  the second Ricci-Chern curvature $ Tr_\omega R^E$ is nonpositive everywhere and negative at some point,  then $E$ has no holomorphic section, i.e. $H^0(M,E)=0$;

\item[(3)]  If  the second Ricci-Chern curvature $ Tr_\omega R^E$ is $p$-nonpositive everywhere and $p$-negative at some point, then $\Lambda^{q}E$ has no holomorphic section for  any $p\leq p\leq
rank(E)$. \ed \ecorollary

Now we can apply it to the tangent and cotangent bundles of compact
Hermitian manifolds.

\bcorollary\label{hermitianvanishing} Let $(M,\omega)$ be a compact
Hermitian manifold and $\Theta$ is the Chern curvature of the Chern
connection $\nabla^{CH}$ on the holomorphic tangent bundle
$T^{1,0}M$. \bd\item[(1)] If  the second Ricci-Chern curvature
$\Theta^{(2)}$ is nonpositive everywhere and negative at some point,
then $M$ has no holomorphic vector field, i.e. $H^0(M,T^{1,0}M)=0$;

\item[(2)] If  the second Ricci-Chern
curvature $\Theta^{(2)}$ is nonnegative everywhere and positive at
some point, then $M$ has no holomorphic $p$-form for any $1\leq
p\leq n$, i.e. $H^{p,0}_{\bp}(M)=0$; In particular, the arithmetic
genus \beq \chi(M,\sO)=\sum(-1)^p h^{p,0}(M)=1 \eeq
\item[(3)] If  the second Ricci-Chern
curvature $\Theta^{(2)}$ is $p$-nonnegative everywhere and
$p$-positive at some point, then $M$ has no holomorphic $q$-form for
any $p\leq q\leq n$, i.e. $H^{q,0}_{\bp}(M)=0$. In particular, if
the scalar curvature $S^{CH}$ is nonnegative everywhere and positive
at some point, then $H^{0}(M, mK_M)=0$ for all $m\geq 1$ where $K_M$
is the canonical line bundle of $M$. \ed

\bproof Let $E=T^{1,0}M$ and $h$ be an Hermitian metric on $E$ such
that the second Ricci-Chern curvature $Tr_{\omega_h}\Theta$ of
$(E,h)$ satisfies the assumption. It is obvious that all section
spaces in consideration are independent of the choice of the metrics
and connections.

The metric on the vector bundle $E$ is fixed.
 Now we choose a Gauduchon metric $\omega_G=e^u\omega_h$ on
$M$. Then the second Ricci-Chern curvature
$\tilde\Theta^{(2)}=Tr_{\omega_G}\Theta=e^{-u}Tr_{\omega_h} \Theta$
shares the semi-definite property with
$\Theta^{(2)}=Tr_{\omega_h}\Theta$. For the safety, we repeat the
 arguments in Theorem \ref{vanishing} briefly. If $s$ is a holomorphic
section of $E$, i.e., $\bp_Es=\bp s=0$, by formula
\ref{gauduchonlaplaciancompare}, we obtain
 \beq 0=\|\p_E s\|^2+\left(\sq\left[\p_E\bp_E+\bp_E\p_E,\Lambda_
G\right]s,s\right)= \|\p_E s\|^2-\left(Tr_{\omega_G}\Theta
s,s\right)\eeq If $Tr_\omega\Theta$ is nonpositive everywhere, then
$\p_E s=0$ and so $\nabla^E s=0$. If $Tr_\omega\Theta$ is
nonpositive everywhere and negative at some point, we get $s=0$,
therefore $H^{0}(M,T^{1,0}M)=0$. The proofs of the other parts are
similar. \eproof

\ecorollary

%

\bremark  It is well-known that the first Ricci-Chern curvature
$\Theta^{(1)}$ represents the first Chern class of $M$. But on an
Hermitian manifold, it is possible that the second Ricci-Chern
curvature $\Theta^{(2)}$ is not in the same $(d,\p,\bp)$-cohomology
class as $\Theta^{(1)}$. For example, $\S^3\times \S^1$ with
canonical metric has strictly positive second Ricci-Chern curvature
but it is well-known that it has vanishing first Chern number
$c_1^2$. For more details see Proposition \ref{hopf}. Therefore,
$\Theta^{(2)}$ in Proposition \ref{hermitianvanishing} can
\textsc{not} be replaced by $\Theta^{(1)}$. It seems to be an
interesting question:  if $(M,\omega)$ is a compact Hermitian
manifold and its  first Ricci-Chern curvature is nonnegative
everywhere and positive at some point,   is the first Betti number
of $M$ zero? In particular, is it K\"ahler in dimension $2$?

\eremark

As special cases of our results, the following results for K\"ahler
manifolds are well-known, and we list them here for the convenience
of the reader. Let $(M,h,\omega)$ be a compact K\"ahler manifold.

\bd \item[(1)] If the Ricci curvature is nonnegative everywhere,
then any holomorphic $(p,0)$ form is parallel;

\item[(2)] If the Ricci curvature is nonnegative everywhere and positive
at some point, then $h^{p,0}=0$ for $p=1,\cdots, n$. In particular,
the arithmetic genus $\chi(M,\sO)=1$ and $b_1(M)=0$;

\item[(3)] If the scalar curvature is nonnegative everywhere and positive
at some point, then $h^{n,0}=0$.
 \item[(A)] If the Ricci curvature is nonpositive everywhere,
then any holomorphic vector field is parallel;

\item[(B)] If the Ricci curvature is nonpositive  everywhere and negative
at some point, there is no holomorphic vector field. \ed

%
%
%
%
%
%
%

\subsection{Vanishing theorems on special Hermitian manifolds}

Let $(M,h,\omega)$ be a compact Hermitian manifold and $\nabla$ be
the Levi-Civita connection.

\blemma\label{gg} Let $(M,\omega)$ be a compact balanced Hermitian
manifold. For any $(p,0)$-form  $\phi$
 on $M$, \bd \item[(1)] If $\phi$ is holomorphic, then
$\p^*\phi=0$;

\item[(2)] If $\nabla'\phi=0$, then $\p\phi=0$.

\ed \bproof For simplicity, we assume $p=1$. For the general case,
the proof is the same. By Lemma \ref{dual}, we know, for any
$(1,0)$-form $\phi=\phi_idz^i$, \beq \p^*\phi=-h^{i\bar j}\frac{\p
\phi_i}{\p\bar z^j} \eeq where we use the balanced condition
$h^{i\bar j}\Gamma_{i\bar j}^s=0$. If $\phi$ is holomorphic, then
$\frac{\p\phi_i}{\p\bar z^j}=0$, hence $\p^*\phi=0$. On the other
hand, \beq \nabla'\phi=\left(\frac{\p\phi_i}{\p z^j}-\Gamma_{ji}^m
\phi_m\right)dz^j\ts dz^i \eeq If $\nabla'\phi=0$, we obtain
\beq\p\phi=\frac{\p \phi_i}{\p z^j}dz^j\wedge
dz^i=\Gamma_{ji}^m\phi_m dz^j\wedge dz^i=0\eeq

\eproof

\elemma

\btheorem\label{a1} Let $(M,\omega)$ be a compact  balanced
Hermitian manifold with Levi-Civita connection $\nabla$. \bd
\item[(1)] If the Hermitian-Ricci curvature $(R_{i\bar j})$ is $p$-nonnegative everywhere, then
any holomorphic $(q,0)$-form ($p\leq q\leq n$) is $\p$-harmonic; in
particular, $ h^{q,0}(M)\leq h^{0,q}(M)$ for any $p\leq q\leq n$;
\item[(2)] If the Hermitian-Ricci curvature $(R_{i\bar j})$ is $p$-nonnegative everywhere and $p$-positive
at some point, $H^{q,0}_{\bp}(M)=0 $ for any $p\leq q\leq n$;
 \ed
\noindent In particular, \bd\item[(3)]if the Hermitian-Ricci
curvature $(R_{i\bar j})$ is nonnegative everywhere and positive at
some point, then $H^{p,0}_{\bp}(M)=0$,  for $p=1,\cdots, n$ and so
the arithmetic genus $\chi(M,\sO)=1$ and $b_1(M)\leq h^{0,1}(M)$.
\item[(4)]if the Hermitian-scalar curvature $S$ is nonnegative everywhere and
positive at some point, then
$$H^{0}(M, mK_M)=0 \qtq{for any} m\geq 1$$
where $K_M=\det T^{*1,0}M$.
 \ed

\bproof At first, we assume $p=1$ for $\bf{(1)}$ and $\bf{(2)}$. Now
we consider $E=T^{*1,0}M$ with the induced metric connection
$\nabla^E=\hat \nabla$ for $h$ (see \ref{metricconnection}). By
formula \ref{balancedlaplaciancompare}, we have \beq
\|\bp_Es\|^2=\|\p_Es\|^2+\sq\left( \left[ R^{E},
\Lambda\right]s,s\right) \eeq where $ R^{E}$ is the $(1,1)$-part
curvature of $E$ with respect to the connection $ \nabla^E$. More
precisely, \beq R^E=\p_E\bp_E+\bp_E\p_E=-\hat R_{i\bar j k}^\ell
dz^i\wedge d\bar z^j\ts \pzl \ts dz^k\eeq since $E$ is the dual
vector bundle of $T^{1,0}M$ and the $(1,1)$-part of the curvature of
$T^{1,0}M$ is
 \beq \hat
R_{i\bar j k}^\ell dz^i\wedge d\bar z^j\ts dz^k \ts \pzl \eeq
 If $s=f_idz^i$ is a holomorphic $1$-form, i.e. \beq \bp
s=\frac{\p f_i}{\p \bar z^j}d\bar z^j\wedge dz^i=0\eeq then \beq
\bp_E s=\left(\frac{\p f_i}{\p \bar z^j} -f_{k}\Gamma_{\bar ji}^k
\right) d\bar z^j\ts dz^i=-f_{k}\Gamma_{\bar ji}^k d\bar z^j\ts
dz^i\eeq Without loss of generality, we assume $h_{i\bar
j}=\delta_{ij}$ at a given point. By Proposition
\ref{curvaturecompare}, the  quantity \beq |\bp_E
s|^2=\sum_{i,j,t,n}f_i\bar f_n \Gamma_{\bar j t \bar
i}\bar{\Gamma_{\bar j t\bar n}}= \sum_{i,n} \left(\hat
R^{(2)}_{n\bar i}-R_{n\bar i}\right)f_i\bar f_n \eeq On the other
hand \beq \sq \left\la \left[R^E,
\Lambda\right]s,s\right\ra=\sum_{i,n}\hat R^{(2)}_{n\bar i}f_i\bar
f_n \eeq That is \beq |\bp_E s|^2- \sq \left\la \left[R^E,
\Lambda\right]s,s\right\ra=-\sum_{i,n} R_{n\bar i}f_i\bar f_n\leq
0\eeq if the Hermitian-Ricci curvature $(R_{n\bar i})$ of
$(M,h,\omega)$ is nonnegative everywhere. Then we get \beq 0\leq
\|\p_Es\|^2=\|\bp_Es\|^2-\sq\left( \left[R^E,
\Lambda\right]s,s\right)\leq 0 \eeq That is $\p_Es=0$. Since $$\p_E
s=\nabla^{'E}s=\hat\nabla's=\nabla's=\left(\frac{\p f_i}{\p
z^j}-f_\ell\Gamma_{ij}^\ell\right)dz^j\ts dz^i$$ we obtain
$\nabla's=0$. By Lemma \ref{gg}, we know $\Delta_\p s=0$. In
summary, we get \beq H^{1,0}_{\bp}(M) \subset H^{1,0}_\p (M)\cong
H^{0,1}_{\bp}(M) \eeq If the Hermitian-Ricci curvature $(R_{n\bar
i})$ is nonnegative everywhere and positive at some point, then
$f_i=0$ for each $i$, that is $s=0$. So we proved
$H^{1,0}_{\bp}(M)=0$. The general cases follow by the same arguments
as Theorem \ref{vanishing} and Theorem \ref{hermitianvanishing}.
\eproof \etheorem

%
%
%

\noindent The dual  of Theorem \ref{a1} is \btheorem\label{a0} Let
$(M,h,\omega)$ be a compact balanced Hermitian manifold.

\bd \item[(1)] If $2\hat R^{(2)}_{i\bar j}-R_{i\bar j}$ is
nonpositive everywhere, then any holomorphic vector field is
$\nabla'$-closed;

\item[(2)] If $2\hat R^{(2)}_{i\bar j}-R_{i\bar j}$ is nonpositive everywhere and negative
at some point, there is no holomorphic vector field. \ed

\bproof Let $E=T^{1,0}M$ and $\hat \nabla$ the induced connection on
it. If $s=f^i\frac{\p}{\p z^i}$ is a holomorphic section, then \beq
\bp_Es= f^i\Gamma_{\bar j i}^{\ell} d\bar z^j\ts \frac{\p}{\p
z^\ell}\eeq Without loss generality, we assume $h_{i\bar
j}=\delta_{ij}$ at a given point.  By Proposition
\ref{curvaturecompare},
 \be |\bp_Es|^2-\sq\left\la \left[\hat
R^{1,1}, \Lambda\right]s,s\right\ra &=& \left(\hat R^{(2)}_{i\bar
j}-R_{i\bar j}\right)f^i\bar f^j+ \hat R_{i\bar j}^{(2)}f^i\bar f^j\\
&=& \left(2\hat R^{(2)}_{i\bar j}-R_{i\bar j}\right)f^i\bar f^j\ee
By formula \ref{Laplaciancompare}, \beq 0\leq
\|\p_Es\|^2=\|\bp_Es\|^2-\sq\left( \left[\hat R^{1,1},
\Lambda\right]s,s\right) \eeq So if
 $2\hat R^{(2)}_{i\bar j}-R_{i\bar j}$ is nonpositive everywhere,
 $\p_Es=\nabla's=0$. If $2\hat R^{(2)}_{i\bar j}-R_{i\bar j}$ is nonpositive everywhere and negative
 at
some point, there is no holomorphic vector field. \eproof \etheorem

\bremark \bd\item[(1)] It is obvious that the second Ricci-Chern
curvature $\Theta^{(2)}_{k\bar \ell}$ and Hermitian-Ricci curvature
$R_{k\bar \ell}$ can not be compared. Therefore, Theorem
\ref{hermitianvanishing} and Theorem \ref{a1} are independent of
each other. For the same reason, Theorem \ref{hermitianvanishing}
and Theorem \ref{a0} are independent.
\item[(2)] For a special case in Theorem \ref{a1}, if the Hermitian-Ricci curvature $R_{k\bar
\ell}$ is nonnegative everywhere and positive at some point, by
Proposition \ref{positivebalanced}, the manifold $(M,\omega)$ is
Moishezon. It is well-known that every  $2$-dimensional
Moishezon/balanced manifold is K\"ahler, but there are many
Moishezon non-K\"ahler manifolds in higher dimension( See
\cite{Mic}). \ed\eremark

\noindent The following result was firstly obtained in \cite{IP}:
\bcorollary Let $(M,\omega)$ be a compact Hermitian manifold with
$\Lambda(\p\bp\omega)=0$. Let $\nabla^B$ be the Bismut connection on
$T^{1,0}M$. \bd
\item[(1)] If the first Ricci-Bismut curvature $B^{(1)}$ is
nonnegative everywhere, then every holomorphic $(p,0)$-form is
parallel with respect to the Chern connection $\nabla^{CH}$;

 \item[(2)] If the first
Ricci-Bismut curvature $B^{(1)}$ is nonnegative everywhere and
positive at some point, then $M$ has no holomorphic $(p,0)$-form for
any $1\leq p\leq n$, i.e. $H^{p,0}_{\bp}(M)=0$; in particular, the
arithmetic genus $\chi(M,\sO)=1 $.

\item[(3)] If the first Ricci-Bismut curvature $B^{(1)}$ is
$p$-nonnegative everywhere and $p$-positive at some point then $M$
has no holomorphic $(q,0)$-form for any $p\leq q\leq n$, i.e.
$H^{q,0}_{\bp}(M)=0$. In particular, if the scalar curvature
$S^{BM}$ of the Bismut connection is nonnegative everywhere and
positive at some point, then $H^{0}_{}(M, mK_M)=0$ for any $m\geq
1$. \ed

\bproof  By Proposition \ref{chernbismut}, if
$\Lambda(\p\bp\omega)=0$, then \beq B^{(1)}\leq  \Theta^{(2)}\eeq
Now we can apply Corollary \ref{hermitianvanishing} to get $(1)$,
$(2)$ and $(3)$. 
 \eproof \ecorollary

\bremark For more vanishing theorems on special Hermitian manifolds,
one can see \cite{AI}, \cite{IP}, \cite{GI1}, \cite{GI2} and
references therein. \eremark

\section{Examples of non-K\"ahler manifolds with nonnegative curvatures}

Let $M=\S^{2n-1}\times \S^1$ be the standard $n$-dimensional ($n\geq
2$) Hopf manifold. It is diffeomorphic to $\C^n- \{0\}/G$ where $G$
is cyclic group generated by the transformation $z\rightarrow
\frac{1}{2}z$. It has an induced complex structure of $\C^n-\{0\}$.
For more details about such manifolds, we refer the reader to
\cite{KN}. On $M$, there is a natural metric
 \beq
h=\sum_{i=1}^n\frac{4}{|z|^2}dz^i\ts d\bar z^i \eeq  The following
identities follow immediately \beq \frac{\p h_{k\bar \ell}}{\p
z^i}=-\frac{4\delta_{k\ell}\bar z^i}{|z|^4},\ \ \ \ \frac{\p
h_{k\bar \ell}}{\p\bar z^j}=-\frac{4 \delta_{k\ell}z^j}{|z|^4} \eeq
and \beq \frac{\p^2 h_{k\bar \ell}}{\p z^i\p\bar
z^j}=-4\delta_{k\ell}\frac{\delta_{i\bar j}|z|^2-2\bar z^i
z^j}{|z|^6}\label{wkt} \eeq

\bexample[Curvatures of  Chern  connection]\label{HE} Direct
computation shows that, the Chen curvature components are \beq
\Theta_{i\bar j k\bar \ell}=-\frac{\p^2 h_{k\bar \ell}}{\p z^i\p\bar
z^j}+h^{p\bar q}\frac{\p h_{{k\bar q}}}{\p z^i}\frac{\p h_{p\bar
\ell}}{\p\bar z^j}=\frac{4\delta_{kl}(\delta_{ij}|z|^2-z^j\bar
z^i)}{|z|^6} \eeq and the first and second Ricci-Chern curvatures
are \beq \Theta^{ (1)}_{k\bar \ell}=\frac{n\left(\delta_{k\ell}
|z|^2-z^\ell\bar z^k\right)}{|z|^4},\ \ \ \ \Theta^{(2)}_{k\bar
\ell}=\frac{(n-1)\delta_{k\ell}}{|z|^2}\eeq It is easy to see that
the eigenvalues of $\Theta^{(1)}$ are \beq \lambda_1=0,
\lambda_2=\cdots=\lambda_n=\frac{n}{|z|^2} \eeq Hence,
$\Theta^{(1)}$ is nonnegative and $2$-positive everywhere.
 \eexample

\bexample[Curvatures of Levi-Civita connection] Similarly, we have
\beq \Gamma_{ik}^\ell =-\frac{\delta_{i\ell}\bar
z^k+\delta_{k\ell}\bar z^i}{2|z|^2},\ \ \ \ \Gamma_{\bar j
k}^\ell=\frac{\delta_{jk}z^\ell-\delta_{k\ell}z^j}{2|z|^2}\eeq and
\beq \frac{\p\Gamma_{ik}^\ell}{\p \bar
z^j}=-\frac{\delta_{{k\ell}}\delta_{ij}+\delta_{i\ell}\delta_{jk}}{2|z|^2}+\frac{\delta_{i\ell}z^j\bar
z^k+\delta_{k\ell}z^j\bar z^i}{2|z|^4} \eeq

\beq \frac{\p \Gamma_{\bar j k}^\ell}{\p
z^i}=\frac{\delta_{jk}\delta_{i\ell}-\delta_{k\ell}\delta_{ij}}{2|z|^2}-\frac{(\delta_{jk}z^\ell-\delta_{k\ell}z^j)\bar
z^i}{2|z|^4} \eeq  The complexified Riemannian curvature components
are \beq R_{i\bar j k}^\ell=-\left(\frac{\p \Gamma^{\ell}_{ik}}{\p
\bar z^j}-\frac{\p \Gamma^{\ell}_{\bar jk}}{\p z^i}+\Gamma_{
ik}^{s}\Gamma^{\ell}_{\bar js}-\Gamma_{ \bar jk}^{s}\Gamma^{\ell}_{
is}-{\Gamma_{\bar j k}^{\bar s}\Gamma_{i\bar
s}^\ell}\right)=\frac{\delta_{i\ell}\delta_{jk}}{2|z|^2}-\frac{\delta_{i\ell}
z^j\bar z^k+\delta_{jk}z^\ell\bar z^i}{4|z|^4} \eeq and \beq
R_{i\bar j k\bar\ell}=
\frac{2\delta_{i\ell}\delta_{jk}}{|z|^4}-\frac{\delta_{i\ell}
z^j\bar z^k+\delta_{jk}z^\ell\bar z^i}{|z|^6},\ \ \ \ R_{k\bar
\ell}=\frac{\delta_{k\ell}|z|^2-z^\ell\bar z^k}{2|z|^4}\eeq
\eexample

\bexample[ Curvatures of  Bismut connection] By definition
\ref{curvaturedefinition} and Lemma \ref{bismutlc}, we obtain \beq
B_{i\bar
jk}^\ell=\frac{\delta_{jk}\delta_{i\ell}-\delta_{k\ell}\delta_{ij}}{|z|^2}+\frac{\delta_{ij}\bar
z^k z^\ell+\delta_{k\ell}\bar z^i z^j-\delta_{i\ell}\bar
z^kz^j-\delta_{jk}\bar z^iz^\ell}{|z|^4} \eeq Two Ricci curvatures
are \beq B^{(1)}_{i\bar j}=B^{(2)}_{i\bar
j}=\frac{(2-n)(\delta_{ij}|z|^2-\bar z^i z^j)}{4|z|^2} \eeq On the
other hand, by  \ref{wkt},  it is easy to see $\p\bp\omega=0$ and
$B^{(1)}=0$  for $n=2$.
 \eexample

\bproposition\label{hopf} Let $M=\S^{2n-1}\times \S^1$ be the
standard $n$-dimensional ($n\geq 2$) Hopf manifold with canonical
metric $h$,

\bd\item[(1)] $(M,h)$ has positive second Ricci-Chern curvature
$\Theta^{(2)}$;

\item[(2)] $(M,h)$ has nonnegative first Ricci-Chern curvature $\Theta^{(1)}$, i.e.,  $c_1(M)\geq
0$. Moreover, \beq \int_Mc_1^n(M)=0 \eeq

\item[(3)] $(M,h)$ is semi-positive in the sense of
Griffiths, i.e. \beq \Theta_{i\bar j k\bar \ell}u^i\bar u^j v^k\bar
v^\ell\geq 0 \eeq for any $u, v\in \C^n$; \item[(4)] $R_{k\bar
\ell}$ is nonnegative  and $2$-positive everywhere;

\item[(5)] $(M,h)$ has nonpositive and $2$-negative first Ricci-Bismut
curvature. In particular, $(\S^3\times \S^1,\omega)$ satisfies
$\p\bp\omega=0$ and has vanishing first Ricci-Bismut curvature
$B^{(1)}$.

\ed \eproposition

Although we know all  Betti numbers of Hopf manifold
$\S^{2n-1}\times \S^1$, $h^{p,0}$ is not so obvious. \bcorollary Let
$(M,h)$ be $n$-dimensional Hopf manifold with $n\geq 2$, \bd

\item[(1)] $h^{p,0}(M)=0$ for $p\geq 1$ and
$\chi(M,\sO)=1$. In particular, $h^{0,1}(M)\geq 1$.

\item[(2)] $\dim_\C H^0(M,mK)=0$ for any $m\geq 1$ where
$K=\det(T^{*1,0}M)$. \ed\ecorollary

\section{A natural geometric  flow on Hermitian manifolds}

As we discussed in the above sections, on  Hermitian manifolds, the
second Ricci curvature tensors of various metric connections are
closely related to the geometry of Hermitian manifolds. A natural
idea is to define a flow by using second Ricci curvature tensors of
various metric connections. We describe it in the following.

Let $(M,h)$ be a compact Hermitian manifold.  Let $\nabla$ be an
\emph{arbitrary metric connection} on the holomorphic tangent bundle
$(E,h)=(T^{1,0}M,h)$. \beq \nabla:E\>>>\Om^{1}(E) \eeq It has two
components $\nabla^{'}$ and $\nabla^{''}$, \beq
\nabla=\nabla^{'}+\nabla^{''}\eeq $\nabla^{'}$ and $\nabla^{''}$
induce two differential operators \beq \p_E:
\Om^{p,q}(E)\>>>\Om^{p+1,q}(E)\eeq
\beq\bp_E:\Om^{p,q}(E)\>>>\Om^{p,q+1}(E) \eeq Let $R^{E}$ be the
$(1,1)$ curvature of the  metric connection $\nabla$. More precisely
$R^E$ is a representation of $\p_E\bp_E+\bp_E\p_E$. It is easy to
see that \beq R^E\in\Gamma(M,\Lambda^{1,1}T^*M\ts End(E)) \eeq and
locally, we can write it as \beq R^{E}=R_{i\bar j A}^B dz^i\wedge
dz^j\ts e^A \ts e_B \eeq
 Here
we set $e_A=\frac{\p}{\p z^A}, e^B=dz^B$ where $A, B=1,\cdots ,n$,
since the geometric meanings of $j$ and $A$ are different. It is
well-known that a metric connection $\nabla$ is determined by its
Christoffel symbols \beq \nabla_{\frac{\p}{\p
z^i}}e_A=\Gamma_{iA}^Be_B,\ \ \ \ \nabla_{\frac{\p}{\p\bar
z^j}}e_A=\Gamma_{\bar j A}^Be_B \eeq In particular, we don't have
notations such as $\Gamma_{A i}^B$. It is obvious that \beq R_{i\bar
j B}^A=-\frac{\p \Gamma_{i A}^B}{\p\bar z^j}+\frac{\p\Gamma_{\bar j
A}^B}{\p z^i}-\Gamma_{iA}^C\Gamma_{\bar j C}^B+\Gamma_{\bar j
A}^C\Gamma_{i C}^B \eeq We set the second Hermitian-Ricci curvature
tensor of $(\nabla, h)$  as \beq R^{(2)}=h^{i\bar j} R_{i\bar j
A\bar B} e^A\ts \bar e^B\in \Gamma(M, E^*\ts \bar E^*) \eeq In
general we can study a new class of flows on Hermitian manifolds
\beq
\begin{cases}\frac{\p h}{\p t}=\sF(h)+\mu
h\\
h(0)=h_0\end{cases} \eeq where $\sF$ can be a linear combination of
the first and the second Hermitian-Ricci curvature tensors of
different metric connections on $(T^{1,0}M, h)$. For examples,
$\sF(h)=-\Theta^{(2)}$, the second Ricci-Chern curvature tensor of
the Chern connection, and $\sF(h)=-\hat R^{(2)}$, the second
Hermitian-Ricci curvature tensor of the complexified Levi-Civita
connection, or the second Ricci curvature of any other Hermitian
connection. Quite interesting is to take $\sF(h)=
s\Theta^{(1)}+(1-s)\Theta^{(2)}$ as the mixed Ricci-Chern curvature,
or $\sF(h)=B^{(2)}-2\hat R^{(2)}$ where $B^{(2)}$ is the second
Ricci curvature of the Bismut connection. More generally, we can set
$\sF(h)$ to be certain suitable functions on the metric $h$. For
example, if $\sF(h)=\left(\Delta_{h} S\right)h$, the above equation
will be the Hermitian Calabi flows.

The following result holds for quite general $\sF(h)$, but here for
simplicity we will only take $\sF(h)=-\Theta^{(2)}$ as an example.

\beq \begin{cases}\frac{\p h}{\p t}=-\Theta^{(2)}+\mu h\\
 h(0)=h_0 \end{cases}\label{HRF2}\eeq
 where $\mu$ is a real parameter.
By formula \ref{2chern}, the second Ricci-Chern curvature tensor has
components \beq \Theta^{(2)}_{k\bar \ell}=h^{i\bar j}\Theta_{i\bar j
k\bar \ell}=-h^{i\bar j}\frac{\p^2 h_{k\bar \ell}}{\p z^i\p\bar
z^j}+h^{i\bar j}h^{p\bar q}\frac{\p h_{k\bar q}}{\p z^i}\frac{\p
h_{p\bar \ell}}{\p\bar z^j} \eeq

\btheorem\label{parabolic} Let $(M,h_0)$ be a compact Hermitian
manifold. \bd\item[(1)] There exists small $\eps$ such that, the
solution of flow \ref{HRF2}  exists for $|t|<\eps$, and it preserves
the Hermitian structure;

\item[(2)] The flow \ref{HRF2} preserves the K\"ahler structure,
i.e., if the initial metric $h_0$ is K\"ahler, then $h(t)$ are also
K\"ahler.\ed \bproof \bf{(1).} Let $\Delta_c$ be the canonical
Laplacian operator on the Hermitian manifold $(M,h)$  defined by
\beq \Delta_c=h^{p\bar q}\frac{\p^2}{\p z^p\p\bar z^q}.\eeq
Therefore, the second Ricci-Chern curvature $-\Theta^{(2)}_{i\bar
j}$ has leading term $\Delta_c h_{i\bar j}$ which is strictly
elliptic. The local existence of the flow \ref{HRF2} follows by
general theory of parabolic PDE, and the solution is an Hermitian
metric on $M$.

\bf{(2)}. The coefficients of the tensor $\p\omega$ are given by
\beq f_{i\bar j k}=\frac{\p h_{i\bar j }}{\p z^k}-\frac{\p h_{k\bar
j}}{\p z^i} \eeq Under the flow \ref{HRF2}, we have
\beq\begin{cases} \frac{\p f_{i\bar j k}}{\p
t}=\frac{\p\Theta^{(2)}_{k\bar j}}{\p
z^i}-\frac{\p\Theta^{(2)}_{i\bar j}}{\p z^k}+\mu f_{i\bar j k}\\
f_{i\bar j k}(0)=0\end{cases} \label{kahlerstructure}\eeq At first,
we observe that $f_{i\bar j k}(t)\equiv 0$ is a solution of
\ref{kahlerstructure}. In fact, if $f_{i\bar j k}(t)\equiv 0$, then
$h_{i\bar j}(t)$ are K\"ahler metrics, and so $$\Theta^{(2)}_{i\bar
j}=\Theta^{(1)}_{i\bar j}=-\frac{\p^2\log \det(h_{m\bar n})}{\p
z^i\p\bar z^j}$$ Therefore, \beq \frac{\p\Theta^{(2)}_{k\bar j}}{\p
z^i}-\frac{\p\Theta^{(2)}_{i\bar j}}{\p z^k}=-\frac{\p^3
\log\det(h_{m\bar n})}{\p z^i\p z^k\p\bar z j}+\frac{\p^3
\log\det(h_{m\bar n})}{\p z^i\p z^k\p\bar z j}=0 \eeq On the other
hand, \beq \frac{\p\Theta^{(2)}_{k\bar j}}{\p
z^i}-\frac{\p\Theta^{(2)}_{i\bar j}}{\p z^k}=\Delta_c\left(f_{i\bar
j k}\right)+\qtq{lower order terms} \eeq Hence the solution of
\ref{kahlerstructure} is unique.
 \eproof \etheorem
\bremark Theorem \ref{parabolic} holds also for quite general
$\sF(h)$ which we will study in detail in a subsequent paper
\cite{LLY}. \eremark

The flow \ref{HRF2} has close connections to  several important
geometric flows:

\begin{enumerate}
\item It is very similar to the Hermitian Yang-Mills flow on holomorphic vector
bundles. More precisely, if the flow \ref{HRF2} has long time
solution and it converges to an Hermitian metric $h_{\infty}$ such
that \beq \Theta^{(2)}_{i\bar j}=\mu h_{i\bar j} \eeq The Hermitian
metric $h_\infty$ is Hermitian-Einstein. So, by \cite{LY}, the
holomorphic tangent bundle $T^{1,0}M$ is stable. As shown in Example
\ref{HE}, the Hopf manifold $\S^{2n+1}\times \S^1$ is stable for any
$n\geq 1$. In fact, in the definition of $\Theta^{(2)}_{i\bar j}$,
if we take trace by using the initial metric $h_0$, then we get the
original Hermitian-Yang-Mills flow equation.

\item If the initial metric is K\"ahler, then this flow is reduced to the
usual K\"ahler-Ricci flow(\cite{Cao}).

\item The flow \ref{HRF2} is similar to the harmonic map flow equation as
shown in Theorem \ref{parabolic}. It is strictly parabolic, and so
the long time existence depends on certain curvature condition of
the target manifold as discussed in the pioneering work of
Eells-Sampson in \cite{ES}. The long time existence of this flow and
other geometric properties of our new flow will be studied in our
subsequent work.
\end{enumerate}

\noindent Certain geometric flows and related results have been
considered on Hermitian manifolds recently, we refer the reader to
\cite{ST1}, \cite{ST2}, \cite{ST3} and \cite{Gill}.

\section{Appendix: The proof of the refined Bochner formulas}

\blemma\label{c} On a compact Hermitian manifold $(M,h,\omega)$, we
have \beq [\Lambda, 2\p\omega]=A+B+C \eeq where \beq \begin{cases}
A=-h^{k\bar \ell}h_{i\bar m}\Gamma_{s\bar
\ell}^{\bar m}dz^s\wedge dz^i\ck \\
\bar A^*=-h^{s\bar t}\Gamma_{s\bar k}^{\bar i} d \bar z^k \cbi\cbt
\end{cases}\eeq

\beq \begin{cases} B=-2\Gamma_{i\bar j}^{\bar \ell}dz^i\wedge d\bar
z^j \cbl\\ \bar B^*=2h^{p\bar j}\Gamma_{\ell\bar j}^{\bar s}dz^\ell
\cp\cbs
\end{cases}\eeq

 \beq \begin{cases} C=\Lambda(2\p\omega)=2\Gamma_{j\bar
\ell}^{\bar \ell}dz^j\\
\bar C^*=2 h^{j\bar \ell}\Gamma_{j\bar s}^{\bar s}\cbl=-2h^{j\bar
i}\Gamma_{j\bar i}^{\bar \ell}\cbl
 \end{cases} \eeq
Moreover,
 \bd\item[(1)] $[\Lambda, A]=-\sq \bar
B^*$;

\item[(2)] $[\Lambda, B]=-\sq (2\bar A^*+\bar B^*+\bar C^*)$;

\item[(3)] $[\Lambda, C]=-\sq\bar C^*$.


\ed \bproof All formulas follow by direct computation. \eproof
\elemma

\bdefinition  With respect to  $\nabla'$ and $\nabla''$, we define
\beq
\begin{cases}
D':=dz^i\wedge \nabla'_i\\
D'':=d\bar z^j\wedge \nabla''_{\bar j}
\end{cases}
\eeq The dual operators of $\p,\bp,D',D''$ with respect to the norm
in \ref{norm} are denoted by $\p^*,\bp^*, \delta',\delta''$ and
define \beq\begin{cases}
\delta_0':=-h^{i\bar j} \ci\nabla''_{\bar j}\\
\delta_0'':=-h^{j\bar i}\cbi\nabla'_{j}
\end{cases}\eeq
where $I$ the contraction operator and $I_i=I_{\frac{\p}{\p z^i}}$
and $I_{\bar i}=I_{\bpzi}$. \edefinition

\bremark It is obvious that these first order differential operators
$D', D'', \delta_0'$ and $\delta_0''$ are well-defined and they
don't depend on the choices of holomorphic frames. If $(M,h)$ is
K\"ahler, $D'=\p$, $D''=\bp$, $\delta_0'=\delta'=\p^*$ and
$\delta_0''=\delta''=\bp^*$. \eremark

\blemma\label{G} In the local holomorphic coordinates, \beq
\p=D'-\frac{B}{2} \qtq{and} \bp=D''-\frac{\bar B}{2}\eeq \bproof We
only have  to check them on functions and $1$-forms.  \eproof
\elemma

\blemma\label{dual} On a compact Hermitian manifold $(M,h)$, we have
\beq\begin{cases}
\delta''=\delta_0''-\frac{\bar C^*}{2} \\
\delta'=\delta_0'-\frac{C^*}{2}\end{cases} \label{D*}\eeq For $\p$
and $\bp$, we have \beq
\begin{cases}
\p^*=\delta_0'-\frac{B^*+C^*}{2}\\
\bp^*=\delta''_0-\frac{\bar B^*+\bar C^*}{2}
\end{cases}\eeq
\bproof For any $\phi\in\Om^{p,q-1}(M)$ and $\psi\in \Om^{p,q}(M)$,
by stokes' theorem \be 0&=&\int_M
\bp(\phi\wedge *\bar\psi)\\
&=&\int_M \frac{\p}{\p \bar z^j}\left(  d\bar z^j\wedge \phi\wedge *\bar \psi\right) \\
&=& \int_M \frac{\p}{\p \bar z^j}\left( \la d\bar z^j\wedge \phi,
\psi\ra\frac{\omega^n}{n!}\right)\\
&=&\int_M \frac{\p}{\p \bar z^j}\left( \left\la \phi, h^{j\bar i}
\cbi\psi\right\ra\frac{\omega^n}{n!}\right)\\
&=&\int_M \left( \left\la \nabla''_{\bar j}\phi, h^{j\bar i}
\cbi\psi\right\ra+ \left\la \phi, \nabla'_jh^{j\bar i}
\cbi\psi\right\ra+\left\la \phi, h^{j\bar i} \cbi\psi\right\ra
\frac{\p \log\det
(h_{m\bar n})}{\p \bar z^j} \right)\frac{\omega^n}{n!}\\
&=&\int_M \left( \left\la d\bar z^j\wedge \nabla''_{\bar j}\phi,
\psi\right\ra+ \left\la \phi,h^{j\bar i} \nabla'_j
\cbi\psi\right\ra+ \left\la \phi,\frac{\p h^{j\bar i}}{\p z^j}
\cbi\psi\right\ra+\left\la \phi, h^{j\bar i} \cbi\psi\right\ra
\frac{\p \log\det (h_{m\bar n})}{\p \bar z^j}
\right)\frac{\omega^n}{n!}\\\ee That is \beq (D''\phi,\psi)=\left(
d\bar z^j\wedge \nabla''_{\bar j}\phi, \psi\right)=-\left(
\phi,h^{j\bar i} \nabla'_j \cbi\psi\right) -\left( \phi,\left(
\frac{\p h^{j\bar i}}{\p z^j}+h^{j\bar i} \frac{\p \log\det
(h_{m\bar n})}{\p z^j}\right) \cbi\psi\right) \label{11}\eeq Now we
will compute the second and third terms on the right hand side. \beq
\frac{\p h^{j\bar i}}{\p z^j}+h^{j\bar i} \frac{\p \log\det
(h_{m\bar n})}{\p z^j}=h^{j\bar i}h^{s\bar t}\left(\frac{\p h_{s\bar
t}}{\p z^j} -\frac{\p h_{j\bar t}}{\p z^s}\right)=2h^{j\bar
i}\Gamma_{j\bar t}^{\bar t}=-2 h^{j\bar
 \ell}\Gamma_{j\bar \ell}^{\bar i} \label{12}\eeq On the other hand
\begin{eqnarray} -h^{j\bar i} \nabla'_j \cbi&=&\nonumber-h^{j\bar i} \cbi\nabla'_j
-h^{j\bar
i} I\left(\nabla'_j\frac{\p}{\p\bar z^i}\right)\\
&=&\delta''_0-h^{j\bar i}\Gamma_{j\bar i}^{\bar
\ell}\cbl\label{13}\end{eqnarray} In summary, by \ref{11}, \ref{12}
and \ref{13}, the adjoint operator $\delta''$ of $D''$ is
$$   \delta''=\left(\delta''_0-h^{j\bar i}\Gamma_{j\bar i}^{\bar \ell}\cbl\right)+2h^{j\bar i}\Gamma_{j\bar i}^{\bar \ell}\cbl=\delta_0''-\frac{\bar
C^*}{2}$$ Since $\bp=D''-\frac{\bar B}{2}$, we get
$$\bp^*=\delta''-\frac{\bar B^*}{2}=\delta_0''-\frac{\bar B^*+\bar C^*}{2}$$
\eproof \elemma

\blemma\label{commutator} On a compact Hermitian manifold $(M,h)$,
we have \beq
\begin{cases}
\left[\Lambda,D'\right]=\sq \left(\delta''+\frac{\bar C^*}{2}\right)\\
\left[\Lambda, D''\right]=-\sq (\delta'+\frac{C^*}{2})
\end{cases}
\qtq{and}
\begin{cases}[\delta'',L]=\sq (D'+\frac{C}{2})\\
[\delta',L]=-\sq(D''+\frac{\bar C}{2})
\end{cases}
\eeq 
\bproof By definition \be (\Lambda D' )\phi&=&\left(\sq h^{i\bar
j}\ci\cbj\right)(dz^k\wedge \nabla'_k\phi)\\&=& -\sq h^{i\bar
j} \ci\left(dz^k\wedge \cbj\nabla'_k\phi\right)\\
&=&-\sq h^{i\bar j}\cbj\nabla'_i\phi+\sq h^{i\bar j}dz^k
\ci\cbj\nabla'_k\phi\\&=&\sq \delta_0''+ dz^k\wedge
\nabla'_k\left(\sq h^{i\bar j}\ci\cbj\phi\right)\\&=& \sq
\delta_0''+D'\Lambda\phi\ee where we use the metric compatible
condition \beq \nabla'\omega=0\Longrightarrow \nabla_k'(\Lambda
\phi)=\Lambda(\nabla_k'\phi)\eeq \eproof \elemma

\noindent \blemma\label{DemaillyA} On a compact Hermitian manifold
$(M,h)$, we have \beq
\begin{cases}
\left[\Lambda,\p\right]=\sq \left(\bp^*+\bar\tau^*\right)\\
\left[\Lambda, \bp\right]=-\sq (\p^*+\tau^*)
\end{cases}
\eeq For the dual case, it is \beq
\begin{cases}[\bp^*,L]=\sq (\p+\tau)\\
[\p^*,L]=-\sq(\bp+\bar\tau)
\end{cases}
\eeq 

\bproof By Lemma \ref{commutator},  \ref{G} and  \ref{c}, \be
[\Lambda,
\p]&=&[\Lambda,D']-\left[\Lambda,\frac{B}{2}\right]\\
&=&\sq\left(\delta_0''+\frac{2\bar A^*+\bar B^*+\bar
C^*}{2}\right)\\&=&\sq\left(\delta''+\frac{\bar C^*}{2}+\frac{2\bar
A^*+\bar B^*+\bar C^*}{2}\right)\\
&=&\sq (\bp^*+\bar\tau^*)
 \ee
 The other relations follow by complex conjugate and adjoint
 operations.
\eproof \elemma

\blemma\label{balancedcondition} On an  Hermitian manifold
$(M,h,\omega)$, \beq \bp^*\omega=\sq \Lambda(\p\omega)=\sq
\Gamma_{\ell \bar j}^{\bar j}dz^\ell\eeq

\bproof We have
$$\frac{C}{2}=\Lambda(\p\omega)=\Gamma_{j\bar
\ell}^{\bar \ell}dz^j$$ On the other hand, by Lemma \ref{dual} and
$\delta_0''\omega=0$ \be \bp^*\omega&=&\left(\delta''_0-\frac{\bar
B^*+\bar C^*}{2}\right)\omega= -\frac{\bar B^*\omega}{2}-\frac{\bar
C^*}{2}\omega\\&=& \left(h_{\ell\bar k}h^{p\bar j}h^{i\bar
s}\Gamma_{i\bar j}^{\bar k}dz^\ell
\cp\cbs\right)\left(\frac{\sq}{2}h_{m\bar n}dz^m\wedge d\bar
z^n\right)-\frac{\bar
C^*}{2}\omega\\
&=&-\frac{\sq}{2} h_{\ell\bar k}h^{i\bar j}\Gamma_{i\bar j}^{\bar
k}d z^\ell-\frac{\bar C^*}{2}\omega\\ &=&\frac{\sq}{2} \Gamma_{\ell
\bar j}^{\bar j}dz^\ell-\frac{\bar
C^*}{2}\omega\\
&=&\sq \Gamma_{\ell \bar j}^{\bar j}dz^\ell\\&=&\sq
\Lambda(\p\omega)\ee  \eproof \elemma

Now we assume $E$ is an Hermitian {complex} vector bundle or a
Riemannian vector bundle over a compact Hermitian manifold
$(M,h,\omega)$ and $\nabla^E$ is a metric connection on $E$.

\blemma
 We have the following formula:
 \beq \bp_E^*(\phi\ts s)=(\bp^*\phi)\ts s-h^{i\bar j}\left(\cbj\phi\right) \wedge \nabla_i^Es \label{bpe*}\eeq
 for any $\phi\in \Om^{p,q}(M)$ and $s\in \Gamma(M,E)$.
\bproof The proof of  is the same  as Lemma \ref{dual}. \eproof
\elemma

 \blemma\label{LYA} If $\tau$ is the operator of type $(1,0)$ defined by
$\tau=[\Lambda,2\partial\omega]$ on $\Om^{\bullet}(M,E)$, then\bd
\item[(1)] $[\bp_E^*,L]=\sq(\p_E+\tau)$;
\item[(2)] $[\p^*_E,L]=-\sq(\bp_E+\bar{\tau})$; \item[(3)]
$[\Lambda,\p_E]=\sqrt{-1}(\bp_E^*+\bar{\tau}^{*})$ ;
\item[(4)]
$[\Lambda,\bp_E]=-\sqrt{-1}(\p_E^*+\tau^{*})$.\ed 
%
%
%

\bproof We only have to prove $\bf{(3)}$. For any $\phi\in
\Om^{\bullet}(M)$ and $s\in \Gamma(M,E)$, \be (\Lambda\p_E)(\phi\ts
s)&=&\Lambda\left(\p\phi\ts s+(-1)^{|\phi|}\phi\wedge
\p_Es\right)\\&=& (\Lambda \p \phi)\ts
s+(-1)^{|\phi|}\sq h^{k\bar \ell}\ck\cbl\left(\phi\wedge \p_Es\right)\\
&=& (\Lambda \p \phi)\ts
s+(-1)^{|\phi|}\sq h^{k\bar \ell}\ck\left(\left(\cbl\phi\right)\wedge \p_Es\right)\\
&=&  (\Lambda \p \phi)\ts
s+(-1)^{|\phi|}\sq h^{k\bar \ell}\left(\ck\left(\cbl\phi\right)\right)\wedge \p_Es-\sq h^{k\bar \ell}\cbl(\phi)\wedge \ck\p_Es\\
 &=&  (\Lambda \p \phi)\ts
s+(-1)^{|\phi|}(\Lambda\phi)\wedge \p_Es-\sq h^{k\bar
\ell}\cbl(\phi)\wedge \nabla_k^{E}s
 \ee On the other hand
\be (\p_E\Lambda )(\phi\ts s)&=&\p_E\left((\Lambda \phi)\ts
s\right)\\
&=&(\p\Lambda \phi)\ts s+(-1)^{|\phi|}(\Lambda\phi)\wedge \p_Es \ee
Therefore \be [\Lambda,\p_E](\phi\ts
s)&=&\left([\Lambda,\p]\phi\right)\ts s -\sq h^{k\bar
\ell}\cbl(\phi)\wedge \nabla_k^Es\\
&=&\sq \left(\left(\bp^*+\bar\tau^*\right)\phi\right)\ts s -\sq
h^{k\bar \ell}\cbl(\phi)\wedge \nabla_k^Es \\
&=& \sq \left(\bp_E^*+\bar\tau^*\right)(\phi\ts s)\ee where the last
step follows by \ref{bpe*}.\eproof \elemma


\end{document}